
\def\LaTeX{%
  \let\Begin\begin
  \let\End\end
  \let\salta\relax
  \let\finqui\relax
  \let\futuro\relax}

\def\UK{\def\our{our}\let\sz s}
\def\USA{\def\our{or}\let\sz z}



\LaTeX

\USA


\salta

\documentclass[twoside,12pt]{article}
\setlength{\textheight}{24cm}
\setlength{\textwidth}{16cm}
\setlength{\oddsidemargin}{2mm}
\setlength{\evensidemargin}{2mm}
\setlength{\topmargin}{-15mm}
\parskip2mm


\usepackage[usenames,dvipsnames]{color}
\usepackage{amsmath}
\usepackage{amsthm}
\usepackage{amssymb,bbm}
\usepackage[mathcal]{euscript}

\usepackage{hyperref}
\usepackage{enumitem}

%
%

%
 
\definecolor{rosso}{rgb}{0.85,0,0}
\definecolor{electricviolet}{rgb}{0.56, 0.0, 1.0}
\definecolor{pinegreen}{rgb}{0.0, 0.47, 0.44}

\def\juerg #1{{\color{red}#1}}
\def\juerg #1{{#1}}
\def\an #1{{\color{rosso}#1}}
\def\an #1{{#1}}
\def\pier #1{{\color{blue}#1}} 
\def\pier #1{{#1}} 
%

%




\bibliographystyle{plain}


%
\newtheorem{theorem}{Theorem}[section]
\newtheorem{remark}[theorem]{Remark}
\newtheorem{corollary}[theorem]{Corollary}

\newtheorem{proposition}[theorem]{Proposition}

\finqui
\let\non\nonumber




\def\step #1 \par{\medskip\noindent{\bf #1.}\quad}
\def\jstep #1: \par {\vspace{2mm}\noindent\underline{\sc #1 :}\par\nobreak\vspace{1mm}\noindent}

\def\Lip{Lip\-schitz}
\def\Holder{H\"older}

\def\lhs{left-hand side}
\def\rhs{right-hand side}


\def\multibold #1{\def\arg{#1}%
  \ifx\arg\pto \let\next\relax
  \else
  \def\next{\expandafter
    \def\csname #1#1#1\endcsname{{\bf #1}}%
    \multibold}%
  \fi \next}

\def\pto{.}

\def\multical #1{\def\arg{#1}%
  \ifx\arg\pto \let\next\relax
  \else
  \def\next{\expandafter
    \def\csname cal#1\endcsname{{\cal #1}}%
    \multical}%
  \fi \next}


\def\multimathop #1 {\def\arg{#1}%
  \ifx\arg\pto \let\next\relax
  \else
  \def\next{\expandafter
    \def\csname #1\endcsname{\mathop{\rm #1}\nolimits}%
    \multimathop}%
  \fi \next}

\multibold
qwertyuiopasdfghjklzxcvbnmQWERTYUIOPASDFGHJKLZXCVBNM.

\multical
QWERTYUIOPASDFGHJKLZXCVBNM.

\multimathop
diag dist div dom mean meas sign supp .


\def\Accorpa #1#2 #3 {\gdef #1{\eqref{#2}--\eqref{#3}}%
  \wlog{}\wlog{\string #1 -> #2 - #3}\wlog{}}


\def\<#1>{\mathopen\langle #1\mathclose\rangle}
\def\norma #1{\mathopen \| #1\mathclose \|}

\def\I2 #1{\int_{Q_t}|{#1}|^2}
\def\IT2 #1{\int_{Q_t^T}|{#1}|^2}
\def\IO2 #1{\norma{{#1(t)}}^2}
\def\IOT2 #1{\norma{{#1(T)}}^2}
\def\ov #1{{\overline{#1}}}

\def\SAL{{\cal S}_\gamma}
\def\SO{{\cal S}_0}

\def\fal{F_{1,\gamma}}
\def\faln{F_{1,\gamma_n}}

\def\iot {\int_0^t}

\def\intQt{\int_{Q_t}}

\def\iO{\int_\Omega}

\def\Qtt{\int_{Q_t^T}}

\def\dt{\partial_t}
\def\ddt{\partial_{tt}}
\def\dn{\partial_{\bf n}}
\def\checkmmode #1{\relax\ifmmode\hbox{#1}\else{#1}\fi}

\def\indi{I_{[-1,1]}}
\def\sindi{\partial I_{[-1,1]}}

\def\erre{{\mathbb{R}}}

\def\enne{{\mathbb{N}}}

\def\Vp{{V^*}}

\def\J{{\cal J}}

\def\S{{\cal S}}

\def\UU{{\cal U}}

\def\CP{{${\bf (CP)}$}}
\def\CPal{{${\bf (CP_{\boldsymbol{\gamma}})}$}}
\def\CP0{{${\bf (CP_0)}$}}
\def\CPaltil{{${\bf (\widetilde{CP}_{\boldsymbol{\gamma}})}$}}

\def\aeQ{\checkmmode{a.e.\ in~$Q$}}




\def\genspazio #1#2#3#4#5{#1^{#2}(#5,#4;#3)}
\def\spazio #1#2#3{\genspazio {#1}{#2}{#3}T0}

\def\L {\spazio L}
\def\H {\spazio H}
\def\W {\spazio W}

\def\C #1#2{C^{#1}([0,T];#2)}



\def\Lx #1{L^{#1}(\Omega)}
\def\Hx #1{H^{#1}(\Omega)}

\def\Ldue{\Lx 2}

\def\Huno{\Hx 1}
\def\Hdue{\Hx 2}




\let\vp\varphi

\def\a{\alpha}	
\def\b{\beta}

\def\th{\theta}

\def\ph{\varphi}

\def\thc{\theta_c}

\let\TeXchi\chi                         
\newbox\chibox
\setbox0 \hbox{\mathsurround0pt $\TeXchi$}
\setbox\chibox \hbox{\raise\dp0 \box 0 }
\def\chi{\copy\chibox}


\def\fal{F_{1,\gamma}}
\def\faln{F_{1,\gamma_n}}

\def\ubar{\overline{u}}

\def\phial{\varphi_\gamma}
\def\sial{w_\gamma}

\def\phialn{\varphi_{\gamma_n}}
\def\sialn{w_{\gamma_n}}                  
\def\pal{p_\gamma}
\def\qal{q_\gamma}
\def\paln{p_{\gamma_n}}
\def\qaln{q_{\gamma_n}}

\def\Uad{{\cal U}_{\rm ad}}

\def\bvp{\overline\varphi}

\def\ds{{\rm d}s}
\def\bu{{\ov u}}   
\def\bw{{\ov w}}

\usepackage{amsmath}
\DeclareFontFamily{U}{mathc}{}
\DeclareFontShape{U}{mathc}{m}{it}%
{<->s*[1.03] mathc10}{}

\DeclareMathAlphabet{\mathscr}{U}{mathc}{m}{it}
\Begin{document}


%
\title{Optimal control of a nonconserved phase field model\\ of Caginalp type with thermal memory
and\\ double obstacle potential}
\author{}
\date{}
\maketitle
\begin{center}	
\vskip-1.5cm
{\large\sc Pierluigi Colli$^{(1)}$}\\
{\normalsize e-mail: {\tt \href{mailto:pierluigi.colli@unipv.it}{pierluigi.colli@unipv.it}}}\\[0.25cm]
{\large\sc Gianni Gilardi$^{(1)}$}\\
{\normalsize e-mail: {\tt \href{mailto:gianni.gilardi@unipv.it}{gianni.gilardi@unipv.it}}}\\[0.25cm]
{\large\sc Andrea Signori$^{(1)}$}\\
{\normalsize e-mail: {\tt \href{mailto:andrea.signori01@unipv.it}{andrea.signori01@unipv.it}}}\\[0.25cm]
{\large\sc J\"urgen Sprekels$^{(2)}$}\\
{\normalsize e-mail: {\tt \href{mailto:juergen.sprekels@wias-berlin.de}{juergen.sprekels@wias-berlin.de}}}\\[.5cm]
$^{(1)}$
{\small Dipartimento di Matematica ``F. Casorati''}\\
{\small Universit\`a di Pavia}\\
{\small via Ferrata 5, I-27100 Pavia, Italy}\\[.3cm] 
$^{(2)}$
{\small Department of Mathematics}\\
{\small Humboldt-Universit\"at zu Berlin}\\
{\small Unter den Linden 6, D-10099 Berlin, Germany}\\[2mm]
{\small and}\\[2mm]
{\small Weierstrass Institute for Applied Analysis and Stochastics}\\
{\small Mohrenstrasse 39, D-10117 Berlin, Germany}\\[10mm]
\end{center}

\begin{center}
\emph{Dedicated to the memory of Professor Gunduz Caginalp}
\end{center}
\Begin{abstract}
\noindent In this paper, we investigate optimal control problems for a nonlinear state system which 
constitutes a version of the Caginalp phase field system modeling nonisothermal phase transitions with a nonconserved order parameter that takes thermal memory into account. The state system, which is a first-order approximation of a thermodynamically consistent system, is inspired by the theories developed by Green and Naghdi. It consists of
two nonlinearly coupled partial differential equations that govern the phase dynamics and the universal balance law for internal energy, written in terms of the phase variable and the so-called thermal displacement, i.e., a primitive with respect to time of temperature. 
We extend recent results obtained for optimal control problems in which the free energy governing the phase transition 
was differentiable (i.e., of regular or logarithmic type) to the nonsmooth case of a double obstacle potential. 
As is well known, in this nondifferentiable case standard methods to establish the existence of appropriate Lagrange multipliers fail. 
\an{This difficulty is overcome \an{utilizing} of the so-called {\it deep quench} approach. Namely,
the double obstacle potential is approximated} by a family of (differentiable) logarithmic ones for which the existence of optimal controls and first-order necessary conditions of optimality in terms of the adjoint state variables and a variational 
inequality are known. By proving appropriate bounds for the adjoint states of the approximating systems, we can 
pass to the limit in the corresponding first-order necessary conditions, thereby  establishing meaningful first-order necessary optimality conditions also for the case of the double obstacle 
potential.  
\vskip3mm
\noindent {\bf \an{Keywords}:}
phase field model, thermal memory, double obstacle potential, optimal control, first-order necessary optimality conditions, adjoint system, deep quench approximation 

\vskip3mm
\noindent {\bf AMS (MOS) Subject Classification:} {
		35K55, 
        35K51, 
		49J20, 
		49K20, 
		49J50  
		}
\End{abstract}
\salta
\pagestyle{myheadings}
\newcommand\testopari{\sc Colli--Gilardi--Signori--Sprekels}
\newcommand\testodispari{\sc Optimal control of a Caginalp system with double obstacle potential
}
\markboth{\testopari}{\testodispari}
\finqui
%
\section{Introduction}
Suppose that $\Omega\subset \erre^3$ \an{(the two-dimensional case can be treated in the same fashion)} is a bounded, open, and connected set having a smooth boundary
$\Gamma=\partial\Omega$ with unit outward normal field \an{$\nnn$ and associated normal derivative} $\dn$, and let, for some \an{given} final time $T>0$,
$$
\an{Q_t:=\Omega\times (0,t)  \,\mbox{ for }\,t\in(0,T), \quad Q:=Q_T,  \,\mbox{ and }\,
\Sigma:=\Gamma\times(0,T).}
$$ 
We consider in this paper the following optimal control problem:

\vspace{3mm}\noindent
{\bf (CP)} \quad Minimize the cost functional
\begin{align}\label{cost}
		\notag
		 \J((\vp, w),u)  &: = \frac {k_1} 2 \norma{\vp - \vp_Q}_{L^2(Q)}^2
		+ \frac {k_2} 2 \norma{\vp(T) - \vp_\Omega}^2_{\an{\Lx2}}
		+ \frac {k_3} 2  \norma{ w - w_Q}_{L^2(Q)}^2 
		\\ & \quad 	\notag
		+ \frac {k_4} 2  \norma{ w(T) - w_\Omega}^2_{\an{\Lx2}}
		+ \frac {k_5} 2 \norma{\dt w - w_Q'}_{L^2(Q)}^2 
		+ \frac {k_6} 2  \norma{\dt w(T) - w'_\Omega}^2_{\an{\Lx2}}
		\\ & \quad 
		+ \ell\, \norma{u}_{L^2(Q)}^2  
\end{align}
subject to the state system 
\begin{alignat}{2}
\label{ss1}
&\dt\ph-\Delta\ph + F_1'(\ph) + \an{\tfrac 2 {\thc}} F_2'(\ph) - \an{\tfrac 1{\thc^2}}{\dt w}\,F_2'(\ph) = 0 \quad && \mbox{\an{in} }\,Q\,,\\
\label{ss2}
&\ddt w -\Delta (\alpha \dt w  + \beta w) + F_2' (\ph) \dt \ph =u \quad && \mbox{\an{in }}\,Q\,,\\
\label{ss3}
&\dn \ph=\dn (\alpha \dt w + \beta w)=0 \quad && \mbox{\an{on} }\,\Sigma,\\
\label{ss4}
&\ph(0)=\ph_0,\quad w(0)=w_0,\quad \dt w(0)=v_0 \,\quad && \mbox{\an{in} }\,\Omega\an{,}
\end{alignat}
\Accorpa\Sys {ss1} {ss4}
and \an{the} control constraint
\begin{align*}
	u \in \Uad,
\end{align*}
with the control space  $\UU := L^\infty(Q) $ and
\begin{align}
	&\Uad :=
		\big\{ 			
			u \in \UU :  u_* \leq u \leq u^* \,\,\aeQ		\big\}.
		\label{Uad}
\end{align}
Above, the symbols \an{$k_1, ..., k_6,$ and $\ell$}  denote nonnegative constants which are not all zero, 
$\vp_Q, w_Q, w_Q'\in L^2(Q)$ \,and $\,\vp_\Omega,w_\Omega,w_\Omega'\in L^2(\Omega)$  denote some prescribed targets,
and $\theta_c>0$ denotes a critical temperature.
As for the set of admissible controls $\Uad$, we 
assume that $u_*, u^*$ with $u_*\,\le\,u^*$ \an{(to be intended pointwise)} are prescribed threshold functions in $ L^\infty(Q)$.  
Notice that the control variable $\,u\,$ has the physical meaning
of a distributed heat source.

The state system \Sys\
constitutes an extension of the phase field model 
for nonisothermal phase transitions with nonconserved  order parameter taking place
in the container $\Omega$ which was 
introduced by G. Caginalp in his seminal paper \cite{Cag}. 
The primary variables of the system are $\ph$, the  order parameter
of the phase transition, and $w$, the so-called {\em thermal displacement} or 
{\em freezing index}. 
The latter is directly connected to the absolute temperature $\th$ of the system through the relation
\begin{align}
	\label{thermal_disp}
	w (\cdot , t)  = w_0 + \iot \th(\cdot, s) \,\ds, \quad t \in[0,T].
\end{align}
In the recent paper \cite{CSS3}, the system \Sys\ was derived from the general principles of thermodynamics, where the specific free energy governing the evolution was (up to some physical constants that here are assumed
to equal unity \an{for simplicity}) of the form
\begin{equation}
\label{free}
 F(\th,\ph)=
 \th (1-\ln (\th/\th_1)) + \th\,F_1(\ph)+F_2(\ph)  +\frac{\th}2\,|\nabla\varphi|^2\,.
\end{equation} 
Here, $\theta_1>0$ is a reference temperature. In this framework, \an{equation} \eqref{ss1} describes the dynamics of the phase evolution, while \an{equation} \eqref{ss2} is the universal balance of internal energy, in which the heat flux is in place of
the standard Fourier law assumed in
the Green--Naghdi form (see\an{, e.g.,} \cite{GN91,GN92,GN93,PG09}) 
\begin{equation}
\label{Green}
 \mathbf q=-\alpha \nabla (\dt w )- \beta \nabla w\,,\quad\mbox{with positive constants \an{$\a$ and $\b$}}\,,
\end{equation}
which models the presence of a thermal memory in the system.

For the nonlinearities driving the phase transformation, we assume that \pier{$F_2$ is differentiable with a globally \Lip\ continuous derivative $F_2'$} 
on $\erre$ (typically, \pier{$F_2$ is} a concave function), while for $F_1$ we consider  the \an{convex} functions   
\begin{align}
\label{logpot} 
&F_{\rm 1,log}(r)=\left\{\begin{array}{ll}
(1+r)\,\ln(1+r)+(1-r)\,\ln(1-r)\quad&\mbox{for $\,r\in (-1,1)$}\\
2\,\pier{\ln 2}\quad&\mbox{for $\,r\in\{-1,1\}$\,,}\\
+\infty\quad&\mbox{for $\,r\not\in[-1,1]$}
\end{array}\right.\\[1mm]
&\indi (r)=\left\{\begin{array}{ll}
0\quad&\mbox{for $\,r\in [-1,1]$}\\
+\infty\quad&\mbox{for $\,r\not\in[-1,1]$}
\end{array}\right. \,.
\label{obspot}
\end{align}
We assume that $\indi+F_2$ is a double-well potential. This is actually the case 
if $F_2(r)=k(1-r^2)$, where $k>0$; the function
$\indi+F_2$ is then referred to as a 
{\em double obstacle} potential. Note also that $F'_{\rm 1,log}(r)$ becomes 
unbounded as $r\searrow -1$ and ${r\nearrow 1}$,
and that in the case of \eqref{obspot} the first equation \eqref{ss1} has to be interpreted
 as a differential inclusion,
where $F_1'(\ph)$ is understood in  the sense of subdifferentials. Namely, \eqref{ss1} has to be written
as
\begin{equation}
\label{ss2new}
\dt\ph-\Delta\ph + \xi + \frac 2 {\thc} F_2'(\ph) - \frac 1{\thc^2}{\dt w}\,F_2'(\ph) = 0\,, \quad \xi\in\sindi(\ph).
\end{equation}
We also notice that the equation \eqref{ss1} is of Allen--Cahn type and is suited for the case of nonconserved order parameters (while the case of a conserved order parameter would require a Cahn--Hilliard structure). 

As far as well-posedness is concerned, the above model was already treated in \cite{CSS3} (see Theorem~\ref{THM:EX:WEAK} 
and Theorem~\ref{THM:EX:STRONG} below).
A discussion of a simpler problem for \an{\eqref{ss1}--\eqref{ss4}} was already given in~\cite{MQ10}. The papers~\cite{CC12,CC13} dealt with well-posedness
issues and asymptotic analyses with respect to the positive coefficients \an{$ \alpha$ and $\beta$} as one of them approaches zero. Other results for this class of systems 
may be found in \cite{CDM, CGM13}. 
About optimal control problems for phase field systems, in particular of Caginalp type, 
we can quote the pioneering work \cite{HJ92}; one may also see the specific sections in the monograph~\cite{Trol}. For other contributions, we mention the article \cite{LS} dedicated to a thermodynamically consistent version of the phase field system described above, and the more recent papers \cite{CGMR1} and \cite{CGS26}, where the interested reader can find a list of related references. 

The optimal control problem {\bf (CP)} has been treated in \cite{CSS3} for the
case of regular and logarithmic nonlinearities $F_1$. For such nonlinearities, differentiability properties of 
the control-to-state mapping, the existence of optimal controls, 
as well as first-order necessary optimality conditions could be established. Actually, in \cite{CSS3} a more 
general cost functional was considered which involved the initial temperature
$v_0$ as a second control variable.
In this paper, we focus on the nondifferentiable case when $F_1=\indi$. While a well-posedness result \an{for system \Sys}
was proved in \cite{CSS3} also for
this case\an{,} in which \eqref{ss1} has to be replaced by the inclusion \eqref{ss2new}\an{,} the corresponding optimal 
control problem has not yet been treated. While the existence of optimal controls is not too difficult to show, the  
derivation of necessary optimality is challenging since standard constraint qualifications to establish the 
existence of suitable Lagrange multipliers are not available. In order to overcome this difficulty, we 
employ the so-called \emph{deep quench approximation} which has proven to be a useful tool in a number 
of optimal control problems for Allen--Cahn and Cahn--Hilliard systems involving double obstacle potentials\an{: see, e.g.,} the papers
\cite{CFGS, CFS, CGSEECT,CGSconv, CGS24, CGS2021,S_DQ}\an{.}

In all of these works, the starting point was that the optimal control problem had been successfully treated (by proving Fr\'echet differentiability of 
the control-to-state operator and establishing
first-order necessary optimality conditions in terms of a variational inequality and the adjoint
state system) 
for the case when in the state system \Sys\ the nonlinearity 
$F_1$ is, for $\gamma>0$, given by 
\begin{equation*}
F_{1,\gamma}:=\gamma\, F_{\rm 1,log}.
\end{equation*}
We obviously have that
\begin{alignat}{2}
\label{monoal}
&0\,\le\,F_{1,{\gamma_1}}(r)\,\le\,F_{1,{\gamma_2}}(r)\quad && \forall\,r\in\erre, \quad\mbox{if }\,0<\gamma_1<\gamma_2, 
\\[1mm]
\label{limF1al}
&\lim_{\gamma\searrow0} F_{1,\gamma}(r)\,=\,\indi(r) \quad && \forall\,r\in \erre.
\end{alignat}
In addition, {we note that} $\,F_{\rm 1,log}'(r)=\ln\left(\frac{1+r}{1-r}\right)$ \,and\, $F_{\rm 1,log}''(r)=\frac 2{1-r^2}>0$\, for
$r\in (-1,1)$, and thus, in particular,
\begin{align*}
&\lim_{\gamma\searrow0}\, F_{1,\gamma}'(r)= {\lim_{\gamma\searrow0}\, \gamma\, F_{1,\rm log}'(r) = 0} \quad\mbox{for }\,\an{r \in (-1,1)},\\
&\lim_{\gamma\searrow 0} \Bigl(\,\lim_{r\searrow -1}\fal'(r)\Bigr)=-\infty, \quad
\lim_{\gamma\searrow0} \Bigl(\,\lim_{{r\nearrow 1}}\fal'(r)\Bigr)=+\infty.
\end{align*} 
We may therefore regard the graphs of the single-valued functions 
\begin{equation*}
F_{1,\gamma}'(r)\,=\, \gamma\, F_{\rm 1,log}'(r), \quad \mbox{for}\quad r\in (-1,1)\quad\mbox{and}\quad \gamma>0,
\end{equation*} 
as approximations to the graph of the multi-valued subdifferential $\sindi$ 
from the interior of $(-1,1)$.

For both $F_1=\indi$ (in which case \an{\eqref{ss1}} has to be replaced by 
the inclusion \eqref{ss2new}) and $F_1=\fal$ (where $\gamma>0$), the well-posedness 
results from \cite{CSS3} yield the existence of a unique solution
$(\vp,w)$ and $(\phial,\sial)$ to the state system
\Sys. It is  natural to expect that $(\phial,\sial)\to (\vp,w)$ as $\gamma\searrow0$ in a suitable topology.
Below (cf.~Theorem~\ref{THM:DQ}), we will show that this is actually true.
Owing to the \an{above} construction \an{and the singularity of $\fal$}, the approximating functions $\,\phial\,$  automatically attain their values in
the domain of $\indi$; that is, we have $\|\phial\|_{L^\infty(Q)}\,\le\,1$ \,for all $\gamma>0$
\an{which shows that the order parameter $\phial$ is limited to range in the physical interval $[-1,1].$}

In the following the optimal control problem {\bf (CP)} will be denoted by \CP0\ if
$F_1=\indi$ and by \CPal\ if $F_1=\fal$\an{, $\gamma \in (0,1]$}. The general strategy is to derive 
uniform (with respect to $\gamma\in (0,1]$) a priori estimates for the state and adjoint state variables of 
an ``adapted'' version of \CPal\ that are sufficiently strong as to permit a passage to the 
limit as $\gamma\searrow0$ in order to derive meaningful first-order necessary optimality conditions also for \CP0. 
It turns out that this strategy succeeds. 

The remainder of the paper is organized as follows: Section 2 is devoted to \an{the mathematical} analysis of the state system \Sys, where we cite results obtained in \cite{CSS3} and
\an{derive} a qualitative estimate for the difference between the solutions for different values of
$\gamma$.
The subsequent Section 3 then brings a discussion of the deep quench approximation and its properties.
In the final Section 4, the first-order necessary optimality conditions for the problem \CP0\ will be derived.

At this point, we fix some notation we are going to employ throughout the paper.
Given a Banach space $X$, we denote by $\norma{\cdot}_X$ the corresponding norm, by $X^*$ its 
dual space, and by $\< \cdot , \cdot >_X$ the related duality pairing between $X^*$ and $X$. 
The standard Lebesgue and Sobolev spaces defined on $\Omega$, for every $1 \leq p \leq \infty$ and $k \geq 0$, are denoted by $L^p(\Omega)$ and $W^{k,p}(\Omega)$, and the associated norms by $\norma{\cdot}_{L^p(\Omega)}=\norma{\cdot}_{p}$ and $\norma{\cdot}_{W^{k,p}(\Omega)}$, respectively. 
For the special case $p = 2$, these spaces become Hilbert spaces, and we denote by $\norma{\cdot}=\norma{\cdot}_2$ the norm of $\Lx2$ and employ the usual notation $H^k(\Omega):= W^{k,2}(\Omega)$. 

For convenience, we also introduce the \an{shorthands}
\begin{align}
  & H := \Ldue \,, \quad  
  V := \Huno\,,   \quad
  W := \{v\in\Hdue: \ \dn v=0 \,\mbox{ on $\,\Gamma$}\}.
  \label{def:HVW}
\end{align}
Besides, for Banach spaces $X$ and $Y$, we introduce the linear space
$X\cap Y$, which becomes a Banach space when equipped with its natural norm $\,\|v\|_{X \cap Y}:=
\|v\|_X + \|v\|_Y$, for $v\in X\cap Y$. 
To conclude, for \an{a} normed \an{space} $\,X\,$ and $\,v\in L^1(0,T;X)$, we \juerg{introduce the convolution products}
\begin{alignat}{2}
	(1 * v)(t) & :=\int_0^t v(s)\,\ds,\quad && \hbox{$t \in[0,T]$},
	\label{intr1}
	\\
	(1 \circledast v)(t) & :=\int_t^T v(s)\,\ds, \quad && \an{\hbox{$t \in[0,T]$}}.
	\label{intr2}	
\end{alignat}

\section{Properties of the state system}
\label{SEC:WP}
\setcounter{equation}{0}

The following \an{structural} assumptions are postulated throughout this paper.
\begin{enumerate}[label={\bf (A\arabic{*})}, ref={\bf (A\arabic{*})}]
\item \label{ass:1:constants}
\,$\a,\b,$ and $\thc$ are \an{fixed} positive constants.
\item \label{ass:2:F2}
\,$F_2\in C^3(\erre)$, and $F_2'$ is a \Lip\ continuous function on $\erre$.
\item \label{ass:3:initials}
\,$\ph_0 \in V$, \,\,$w_0 \in V, \,\,v_0 \in H$,\,\, and there are constants $r_*,r^*$  such that\\
\hspace*{1mm} $-1<r_*\le \ph_0(x)\le r^*<1$ for almost every $x\in\Omega$. 
\end{enumerate}
The first result concerns the existence of weak solutions. 
\an{Here, let us incidentally notice that conditions \ref{ass:2:F2} and \ref{ass:3:initials} may \juerg{in fact be} weakened if we are merely interested in the existence of weak solutions to system \Sys.
However, as will be clarified later on, the optimal control problem we aim at solving requires sufficient regularity properties for the state system to derive the corresponding first-order necessary conditions.
For this reason, we immediately assume \ref{ass:2:F2} and \ref{ass:3:initials} to hold in the current form.
} 
\begin{theorem}
\label{THM:EX:WEAK}
Assume that \ref{ass:1:constants}--\ref{ass:3:initials} hold, and assume that either $F_1=F_{1,\gamma}$ for some $\gamma\in (0,1]$ 
or $F_1=\indi$. Then the state system \Sys\ has for every $u\in L^2(0,T;H)$ a unique weak solution $(\ph, w, \xi)$ in the sense that 
\begin{align*}
	&\ph \in \H1 H \cap \L\infty V \cap \L2 W,
	\\ 
	& w \in \H2 {\Vp}\cap\W{1,\infty} H \cap \H1 V,
	\\ 
	&\pier{\xi  \in	\L2 H,
	\quad \text{and} \quad \xi \in \partial F_1(\ph) \,\, a.e. \,\, \text{in} \,\, Q,}
\end{align*}
and that the variational equalities
\begin{align}
	\label{wf:1}
	& \iO \dt\ph \, v 
	+ \iO \nabla \ph \cdot \nabla v
	+\iO  \xi v
	+\frac 2 {\thc}\iO  F_2'(\ph) v
	- \frac 1{\thc^2} \iO {\dt w}\,F_2'(\ph)v = 0\,,
	\\ 
	\label{wf:2}
	& \<\ddt w, v>_{V}
	+\a \iO \nabla (\dt w)\cdot \nabla v
	+ \b \iO \nabla  w  \cdot \nabla v
	+ \iO F_2' (\ph) \dt \ph \, v
	=
	\iO {u}  v\,,
\end{align}
are satisfied for every $v \in V$ and almost every $t\in(0,T)$.
Moreover, it holds that
\begin{align*}
	\ph(0)=\ph_0, 	\quad w(0)=w_0, \quad \dt w(0)={v_0}.
\end{align*}
Furthermore, there exists a constant $K_1>0$, which depends only on $\Omega,T,\a,\b,\thc$ and the data of the system, such that
\begin{align}
	\non
	& \norma{\ph}_{\H1 H \cap \L\infty V \cap \L2 {\Hx2}}
	+ \an{\norma{F_1(\ph)}_{\L\infty {\Lx1}}}
	\\
	& \quad \label{weaksol:estimate}
	+ \norma{w}_{\H2{\Vp} \cap\W{1,\infty} H \cap \H1 V}
		\leq K_1\,.
\end{align}
\end{theorem}
\begin{remark}
{\rm In the single-valued situation, when $F=F_{1,\gamma}$, the inclusion $\xi\in\partial F_{1,\gamma}(\vp)$ becomes the
identity $\xi=F_{1,\gamma}'(\vp)$. Note also that the initial conditions are meaningful at least in $H$, since, in particular, $\ph \in \C0V $ and $w\in \C1H$ by interpolation. \pier{Besides, we point out that \eqref{weaksol:estimate} implies that $-1\leq \ph\leq 1$ a.e. in $Q$, which entails that $\ph$ is also uniformly bounded in $L^\infty (Q)$.}}
\end{remark}

\noindent{\em Proof of Theorem \ref{THM:EX:WEAK}} \,\,\,
The existence and uniqueness results follow as a special case of Theorem 2.1 \an{and Theorem 2.2} in \cite{CSS3}. Moreover, it follows from 
\ref{ass:3:initials} that $\indi(\vp_0)=0$ and that, for all $\gamma\in (0,1]$, we have $\|F_{1,\gamma}(\vp_0)\|_
{\an{1}}=\gamma \|F_{1,{\rm log}}(\vp_0)\|_{\an{1}}
\,\le\,2\,\ln(2)\,|\Omega|$, where $\an{|\Omega|}$ denotes the \an{Lebesgue} measure of $\Omega$. Therefore, the estimates 
performed in the proof of \pier{\cite[Thm.~2.1]{CSS3}} apply for both $F_1=\indi$ and $F_1=F_{1,\gamma}$, $\an{\gamma \in (0,1]}$, which proves the 
validity of \eqref{weaksol:estimate}.\qed

\vspace{2mm}
 By virtue of Theorem~\ref{THM:EX:WEAK}, the control-to-state operators $\S_0:u\mapsto (\vp,w,\xi)$ and 
$\S_\gamma:u\mapsto (\vp,w\pier{,F'_{1,\gamma}(\ph)})$  corresponding to the choices $F_1=\indi$
and $F_1=F_{1,\gamma}$ for $\gamma>0$ , respectively, are well defined. 
\an{Obviously, the above result tacitly \pier{defines} the space where the solution operators $\S$ and $\S_\gamma$ map.
Namely, it shows that both the \pier{operators} $\S$ and $\S_\gamma$ have to be intended as \pier{mappings} from the control space $\cal U$ into the regularity space $\cal Y$ defined by
\begin{align*}
	{\cal Y} & := \big( \H1 H \cap \L\infty V \cap \L2 W \big)  
	\\ & \quad
	\times \big( \H2 {\Vp}\cap\W{1,\infty} H \cap \H1 V\big) \pier{{}\times  \L2 H.}
\end{align*}
We will see in the forthcoming Theorem~\ref{THM:EX:STRONG} that, upon requiring stronger \pier{assumptions} for the 
initial data (cf. \eqref{initials}), the solution \juerg{operators} $\S_\gamma$ can also be interpreted as \juerg{mappings} from $\cal U$ into a smaller space than $\cal Y$ defined by the regularity properties \eqref{regphig}--\eqref{regwgt}.}
The following result provides a qualitative comparison between the solutions associated with $F_1=F_{1,\gamma}$
for different values of $\gamma\in (0,1]$.

\begin{theorem}\label{THM:CD:LOG}
Suppose that \ref{ass:1:constants}--\ref{ass:3:initials} hold, and let, for fixed $0<\gamma_1<\gamma_2\le 1$, controls
$u_{\gamma_1},u_{\gamma_2}\in\Uad$ be given. Then the corresponding solutions $(\vp_{\gamma_i},
w_{\gamma_i}  \pier{,F'_{1,\gamma_i}(\vp_{\gamma_i})})=\S_{\gamma_i}
(u_{\gamma_i})$ to {\rm \Sys} associated with $F_1=F_{1,\gamma_i}$ and $u=u_{\gamma_i}$, $i=1,2$, 
satisfy the estimate 
\begin{align}
	& \non
	\|\vp_{\gamma_1}- \vp_{\gamma_2}\|_{\L\infty H \cap \L2 V}
	+ 
	\|w_{\gamma_1}- w_{\gamma_2}\|_{\H1 H \cap \L\infty V	}
	\\ 
	& 
	\quad 
	\leq K_2 \big (\gamma_2-\gamma_1\big)^{1/2} 
	+ K_2 \,\|1*(u_{\gamma_1}-u_{\gamma_2})\|_{\L2 H}\,,	\label{estiga12}
\end{align}
with a positive constant $K_2$ that depends only on 
$\Omega,T,\a,\b,\thc$ and the data of the system.
\end{theorem}

\begin{proof}
We set, for convenience,
\begin{align}
	\label{not:diff:1}
	u&:=u_{\gamma_1}-u_{\gamma_2}, \quad \ph  : = \ph_{\gamma_1} - \ph_{\gamma_2}, 
	\quad 
	w : = w_{\gamma_1} - w_{\gamma_2}, \nonumber\\
	\rho_i &:= F_2'(\ph_{\gamma_i}) \quad \text{for } \, i=1,2,
	\quad	
	\rho  := \rho_1- \rho_2\,.
\end{align}
Using this notation, we take the difference of the weak formulation \eqref{wf:1}--\eqref{wf:2} written for 
$(\ph_{\gamma_i}, w_{\gamma_i}, \xi_{\gamma_i})$, where $\xi_{\gamma_i}=F_{1,\gamma_i}'(\vp_{\gamma_i})$,
$i=1,2$, which yields the system 
\begin{align}
	\label{wf:cd:1}
	& \iO \dt\ph \, v 
	+ \iO \nabla \ph \cdot \nabla v
	+\iO  \bigl(F_{1,\gamma_1}'(\vp_{\gamma_1})-F_{1,\gamma_2}'(\vp_{\gamma_2})\bigr) v
	+\frac 2 {\thc}\iO  \rho \, v\nonumber\\
	&\quad =\, \frac 1{\thc^2} \iO{\dt w}\,\rho_1 v 
	+ \frac 1{\thc^2} \iO {\dt w_2}\,\rho\, v \,, 
	\\ 
	\non
	& \<\ddt w, v>_{V}
	+\a \iO \nabla (\dt w)\cdot \nabla v
	+ \b \iO \nabla  w  \cdot \nabla v
	+ \iO \dt (F_2(\ph_{\gamma_1})-F_2(\ph_{\gamma_2})) v
	\\ & \quad 
	=
	\iO u  v\,,	
	\label{wf:cd:2}
\end{align}
for all $v \in V$ and almost everywhere in $(0,T)$. 
Of course, we also have the initial conditions 
\begin{align}
\label{pier4}
	\ph(0)=0,\quad w(0)=0,\quad \dt w(0)=0 \quad \mbox{a.e. in }\,\Omega.
\end{align}
We now first add the term $\int_\Omega \ph \, v $ to both sides of \eqref{wf:cd:1}, then take 
$v= \ph$  and integrate with respect to time, which leads to the identity
\begin{align}
	& \frac {1}2 \IO2 {\ph}
	+ \int_0^t \norma{\ph(s)}_V^2 \, \ds
	+ \intQt \bigl(F_{1,\gamma_1}'(\vp_{\gamma_1})-F_{1,\gamma_1}'(\vp_{\gamma_2})\bigr)\,\ph \non
	\\ 
	& 	\quad =\,
	 \intQt  \Big( \ph  - \frac{2}{\thc} \rho \Big) \ph
	+ \frac {1}{\thc^2}\intQt {\dt w}\,\rho_1 \, \ph
	+ \frac {1}{\thc^2}\intQt {\dt w_2}\,\rho \,\ph \nonumber\\
	&\qquad -\intQt\bigl(F_{1,\gamma_1}'(\vp_{\gamma_2})-F_{1,\gamma_2}'(\vp_{\gamma_2})\bigr) \vp\,,
	\label{pier3-1}
\end{align}
for all $t\in [0,T]$.
Due to the monotonicity of $F_{1,\gamma_1}'$, we immediately conclude that the third term on the \lhs\ is nonnegative.
The integrals on the \rhs\ have to be estimated individually, where in the following $C>0$ denotes generic
constants that may depend on the data of the system but not on \an{$\gamma_1$ and $\gamma_2$}. 

At first, using the \Lip\ continuity of $\,F_2'\,$ along with the fact that by Theorem~\ref{THM:EX:WEAK} \,$\dt w_i \in \L\infty H \cap \L2 V$, $\ph_i \in \H1 H \cap \L\infty V$, $i=1,2$,  we infer that
\begin{align*}
&{} \intQt  \Big( \ph  - \frac{2}{\thc} \rho \Big) \ph
 \leq C \I2 \ph\,,
 \end{align*} 
 and, with the help of H\"older's inequality and \an{the} continuous embedding $ V \subset L^4(\Omega)$,
\begin{align*} 
&\frac {1}{\thc^2}\intQt {\dt w}\,\rho_1\, \ph
 \,\leq\,  C\int_0^t \norma{\dt w} \Big( \norma{\ph_1}_4 +1\Big) \,\norma{\ph }_4\,{\rm d}s\\ 
&\quad
\leq \, C \Big( \norma{\ph_1}_{\L\infty V} +1\Big)\int_0^t \norma{\dt w} \, \norma{\ph }_V \,{\rm d}s
  \,\leq\,  \frac 1 4 \int_0^t\norma{\ph }_V^2 \,{\rm d}s\,+\, D_1 \int_{Q_t}  |\dt w|^2 \,, 
\end{align*}  
where $D_1$ is a computable and by now fixed constant.
Moreover, we have that 
\begin{align*}  
&\frac {1}{\thc^2}\intQt {\dt w_2}\,\rho \,\ph	 
\,\leq\, C\int_0^t \norma{\dt w_2}_4 \, \norma{\ph} \,\norma{\ph }_4 \,{\rm d}s
	\\  
&\quad\leq\,	C \int_0^t \norma{\dt w_2}_V \, \norma{\ph } \, \norma{\ph }_V	\,{\rm d}s
\,\leq\,  \frac 1 4 \int_0^t\norma{\ph }_V^2 \,{\rm d}s\,+\,C  \int_0^t \norma{\dt w_2}_V^2 \, \norma{\ph }^2 
\,{\rm d}s\,,	
\end{align*}
where, owing to Theorem~\ref{THM:EX:WEAK}, the function $\, t\mapsto \norma{\dt w_2(t) }_V^2\,$ belongs to $L^1(0,T)$. 

Finally, we recall that $F_{1,\gamma_i}=\gamma_i\,F_{1,{\rm log}}$, $i=1,2$, and it follows from the 
convexity of $F_{1,{\rm log}}$ on $[-1,1]$ and the fact that $\,\vp_{\gamma_1}\,,\,\vp_{\gamma_2}\,$ attain
their values in $[-1,1]$ almost everywhere in $Q$ that, a.e. in $Q$,
\begin{align*}
&-\bigl(F_{1,\gamma_1}'(\vp_{\gamma_2})-F_{1,\gamma_2}'(\vp_{\gamma_2})\bigr)\, \vp
\,=\, (\gamma_2-\gamma_1)\,F_{1,{\rm log}}'(\vp_{\gamma_2})\,\vp\,\\
& \quad \le\,
(\gamma_2-\gamma_1)\bigl(F_{1,{\rm log}}(\vp_{\gamma_1})-F_{1,{\rm log}}(\vp_{\gamma_2})\bigr)
\,\le\,(\gamma_2-\gamma_1)\,2\,\pier{{}\ln 2}\,.
\end{align*}
 
Therefore, collecting the above estimates,  it follows from \eqref{pier3-1} that 
\begin{align}
	& \frac {1}2 \IO2 {\ph}
	+ \frac {1}2\int_0^t \norma{\ph(s)}_V^2 \, \ds
	\non \\ 	
	& \quad 
	\leq
	C\,(\gamma_2-\gamma_1) \,
+\,C  \int_0^t \Big(1 +\norma{\dt w_2}_V^2 \Big)\, \norma{\ph }^2\,\ds 
+ D_1 \int_{Q_t}  |\dt w|^2	\,.
	\label{pier3-2}
\end{align}

Next, we integrate  \eqref{wf:cd:2} with respect to time using \eqref{pier4}, then take $v= \dt w$, and integrate once more over $(0,t)$ for an arbitrary $t \in [0,T]$. Adding to both sides the terms 
$\frac {\a}2\, \norma{w(t)}^2 = \a \intQt w \, \dt w$ \an{(recall that now $w(0)=0$ from \eqref{pier4})}, we obtain that            
\begin{align}
	& \I2 {\dt w }
	+ \frac\a2 \norma{w(t)}^2_V
	\,=\,- \b \intQt (1 * \nabla w) \cdot \nabla(\dt  w)
	\non \\ & \quad 
	- \intQt (F_2(\ph_{\gamma_1})-F_2(\ph_{\gamma_2})) \dt w
	+ \intQt (1*u ) \dt w
	+\a \intQt w \, \dt w.
	\label{pier3-3}
\end{align}
We estimate each term on the \rhs\ individually. At first, using the identity
\begin{align*}
	\intQt (1*\nabla w) \cdot \nabla (\dt w )
	= \iO (1*\nabla w(t)) \cdot \nabla w(t)
	- \intQt  |\nabla w|^2, 
\end{align*}
the fact that 
\,$\Vert 1*\nabla w(t) \Vert^2 \leq \Big( \int_0^t \Vert \nabla w 
\Vert \an{\,\ds}\Big)^2\leq T\intQt |\nabla w|^2$, 
as well as Young's inequality, we infer that
\begin{align*}
	- \b \intQt (1 * \nabla w) \cdot \nabla (\dt w )
	\,\leq\, \frac\a 4\, \IO2 {\nabla w}
	+ C \I2 {\nabla w}.
\end{align*}
\an{Then}, we recall that the mean value theorem and the Lipschitz continuity of $F_2'$ yield 
the existence of some $\widehat C>0$ such
that
\begin{align}\label{pier5}
	\left| F_2(r)-F_2(s) \right| 
	\,\leq\,\widehat C (|r| + |s| + 1) |r-s| \quad \hbox{for all } \, r,\, s \in \erre.
\end{align}  
Hence, by virtue of the continuous and compact embedding $V \subset \Lx{p}$, $1\leq p <6 $,
we deduce from Ehrling's lemma (see, e.g., \cite[Lemme~5.1, p.~58]{Lions}), along with 
the H\"older and Young inequalities, and invoking \eqref{weaksol:estimate}, that the second term
on the \rhs\ can be estimated as follows:
\begin{align*}
	&- \intQt \bigl(F_2(\ph_{\gamma_1})-F_2(\ph_{\gamma_2})\bigr) \dt w
	\leq 
	C \iot \big\| |\ph_{\gamma_1}| + |\ph_{\gamma_2}| + 1 \big\|_4 \,
	\norma{\ph_{\gamma_1}- \ph_{\gamma_2}}_4\, \norma{\dt w} \,\ds
	\\ &  \quad 
	\leq 
	\frac 14 \I2 {\dt w }
	+ {C}\bigl(\norma{\ph_{\gamma_1}}_{\L\infty{V}}^2
	\,+\,\norma{\ph_{\gamma_2}}_{\L\infty{V}}^2+1 \bigr) 
	\int_0^t  \|\ph\|_4^2 \,\ds
	\\ &  \quad 
	\leq 
	\frac 14 \I2 {\dt w } + \delta \int_0^t  \|\ph\|_V^2	\juerg{\,\ds}+ C_\delta \I2 \ph	\,,
\end{align*}
for any positive coefficient $\delta$ (yet to be chosen). 
Finally, Young's inequality easily yields that
\begin{align*}
	\intQt (1*u ) \dt w
	+	\a \intQt w \, \dt w 
	& 
	 \leq \frac 14  \I2 {\dt w}
	 + C\I2 {1*{u}}
	+ C \I2 w.
\end{align*}
Thus, in view of \eqref{pier3-3}, upon collecting the above estimates, we realize that 
\begin{align}
	& \frac 12 \I2 {\an{\dt w }}
	+ \frac\a 4\, \norma{w(t)}^2_V
		\non \\ 
		& \quad 
		\leq 	\, 	 
	 \delta \int_0^t  \|\ph\|_V^2 \,\ds	+ C_\delta \I2 \ph	
	 + C\I2 {1*{u}}
	+ C \int_0^t \norma{w}_V^2 \,\ds\,.	
	\label{pier3-4}
\end{align}
At this point, we multiply \eqref{pier3-4} by ${4}D_1 $ and add it to \eqref{pier3-2};
then, fixing $\delta>0 $ such that $4\,D_1\, \delta <1/2$, and applying the Gronwall lemma, we obtain the estimate 
\begin{align*}
	& \norma{\ph}_{\L\infty H \cap \L2 V}
	 + \norma{w}_{\H1 H \cap \L\infty V}
	 \\ & \quad
	 \leq C \bigl((\gamma_2-\gamma_1)^{1/2} \,+\, \norma{1*u}_{\L2 H}\bigr)\,,
\end{align*}
which finishes the proof of the assertion.
\end{proof}

We now derive better regularity and boundedness results for the weak solutions to the state system that correspond 
to the logarithmic potentials $F_1=F_{1,\gamma}$\an{, $\gamma \in (0,1]$}.
\begin{theorem}
\label{THM:EX:STRONG}
Assume that \ref{ass:1:constants}--\ref{ass:3:initials} are fulfilled. Moreover, assume that the following condition is satisfied:
\begin{align}\label{initials}
&\ph_0 \in W,\quad   v_0 \in V\cap L^\infty(\an{\Omega}), \quad w_0\in V\cap L^\infty(\an{\Omega}). 
\end{align}
Then the state system {\rm \Sys} with $F_1=F_{1,\gamma}$ has for every $\gamma\in (0,1]$ and every $u\in {\cal U}$ 
a unique strong solution $(\ph_\gamma, w_\gamma)$ with the regularity 
\begin{align}
	\label{regphig}
	&\ph_\gamma \in {\W{1,\infty} H \cap \H1 V} \cap \L\infty{W},
	\\ 
	\label{regwg}
	&w_\gamma \in {\H2 H \cap \W{1,\infty} V \cap \H1 {W}},
	\\
	\label{regwgt}
	&\partial_t w_\gamma \in L^\infty(Q),
\end{align}
such that the equations \Sys\ are fulfilled almost everywhere in $Q$, on $\Sigma$, or in $\Omega$, respectively.
In addition, 
the phase variable $\ph_\gamma$ enjoys the so-called separation property, i.e., there exist two 
values $r_-(\gamma)\in (-1,r_*\pier{{}]{}}$, $r_+(\gamma)\in \pier{{}[{}}r^*,1)$, which depend only 
on $\Omega,T,\a,\b, \thc$ and the data of the system, such that
\begin{align}
	\label{separation}
	-1<r_-(\gamma) \leq \ph_\gamma \leq  r_+(\gamma) < 1 \quad \hbox{a.e. in } Q.
\end{align}
Moreover, there exists a constant $K_3>0$ such that, for all $\gamma\in (0,1]$,
\begin{align}\label{ssbound}
	\non
	& \norma{\ph_\gamma}_{{\W{1,\infty} H \cap \H1 V} \cap \L\infty{\Hx2}}
	\,+\,\|\fal'(\ph_\gamma)\|_{L^\infty(0,T;H)}
		\\ & 
		\quad 
	+ \norma{w_\gamma}_{\H2 H \cap \W{1,\infty} V \cap \H1 {\Hx2}}
	+ \norma{\dt w_\gamma}_{L^\infty(Q)}
	\leq K_3 \,.
\end{align}
\end{theorem}
\begin{proof}
We want to apply \cite[Thm.~2.3]{CSS3}. To this end, we observe the following facts: first, it obviously holds that 
$\,\Delta \ph_0 - F_{1,\gamma}'(\ph_0) - \tfrac 2{\thc} F_2'(\ph_0) + \tfrac 1 {\thc^2} v_0 
	F_2'(\ph_0)   \in H$. Moreover, the restriction of $F_{1,\gamma}$ to the interval $(-1,1)$ belongs to
$C^2(-1,1)$ and satisfies the conditions $\,\lim_{r\searrow -1}F_{1,\gamma}'(r)=-\infty\,$ and
$\,\lim_{r\nearrow 1}F_{1,\gamma}'(r)=+\infty$.	Therefore, all of the prerequisites for arguing along the 
lines of the proof of \cite[Thm.~2.3]{CSS3} are fulfilled, from which we conclude the validity of the assertion.
We only remark that the global estimate \eqref{ssbound} is a consequence of the special form of 
$F_{1,\gamma}=\gamma F_{1,{\rm log}}$ and of the fact that we only admit
parameters $\gamma$ in the bounded interval $(0,1]$. 
\end{proof}
\begin{remark}
{\rm It cannot be excluded that for $\gamma\searrow0$ we have $\,r_-(\gamma)\to -1\,$ and/or 
$\,r_+(\gamma)\to +1$. Therefore, a global (for $\gamma\in (0,1]$) bound for the $L^\infty(Q)$-norm of $\fal'(\phial)$ cannot be
guaranteed.}
\end{remark}		

\section{Deep quench approximation of the state system}
\setcounter{equation}{0}

In this section, we discuss the deep quench approximation
of the state system \Sys, where we generally assume
that the conditions \ref{ass:1:constants}--\ref{ass:3:initials} and \eqref{initials} are fulfilled. As in the previous section,
we consider the state system for the cases $F_1=\indi$ and  $F_1=\fal$ ($\gamma\in (0,1]$), respectively. 
By Theorem~\ref{THM:EX:STRONG}, we have for every $u\in{\cal U}$ and 
$F_1=\fal$, $\gamma\in (0,1]$, a unique strong solution $(\phial,\sial, \pier{\xi_\gamma})=\SAL(u)$ with
the regularity specified by \eqref{regphig}--\eqref{regwgt} and with 
$\pier{\xi_\gamma :=\fal'(\phial)\in \L \infty H}$, while Theorem~\ref{THM:EX:WEAK} 
implies the existence of a unique weak solution $(\ph^0, w^0, \xi^0)=\SO(u)$ to the weak form 
\eqref{wf:1}--\eqref{wf:2} of the state system for $F_1=\indi$ that enjoys the regularity specified \an{in} Theorem~\ref{THM:EX:WEAK}.  
Clearly, we must have
\begin{equation}
\label{interv}
-1\le \phial\le 1\,\mbox{ a.e. in \,$Q$,} \, \hbox{ for all $\gamma\in (0,1],$ \, and }
\, -1\le \vp^0\le 1\, \mbox{ a.e. in $\,Q$.}
\end{equation} 

We are now going to investigate the behavior of the family
$\{(\phial,\sial)\}_{\gamma>0}$ of deep quench approximations as
$\gamma\searrow0$. We expect that the solution operator $\SAL$ yields an approximation of $\S_0$ as 
$\gamma \searrow 0.$ This is made rigorous through
the following result.
\begin{theorem} \label{THM:DQ}
Suppose that the assumptions \ref{ass:1:constants}--\ref{ass:3:initials} and \eqref{initials}  are fulfilled, and 
let sequences 
$\{\gamma_n\}\subset (0,1]$ and 
$\{u_{\gamma_n}\}\subset \Uad$ be given such that $\gamma_n\searrow0$ and $u_{\gamma_n}
\to  u$ weakly-star in ${\cal U}$  as 
$n\to\infty$ 
for some $ u\in\Uad$. Moreover, let $(\phialn,\sialn\pier{,\xi_{\gamma_n}})=
{\cal S}_{\gamma_n}(u_n)$, $n\in\enne$, and $(\vp^0,\an{w}^0,\xi^0)=\SO(u)$. 
Then, as $n\to\infty$, \pier{we have that}
\begin{align}
\label{conphi}
\phialn\to \varphi^0&\quad\mbox{weakly-star in \,${\W{1,\infty} H \cap \H1 V} \cap \L\infty{W}$}\non\\
&\quad\mbox{and strongly in }
\,C^0(\overline Q), \\
\label{conxi}
\pier{\xi_{\gamma_n}}:= \faln'(\phialn)\to\xi^0&\quad\mbox{weakly-star in }\,L^\infty(0,T;H),\\
\label{conw}
\sialn \to w^0&\quad\mbox{weakly-star in }\,{\H2 H \cap \W{1,\infty} V \cap \H1 {W}}\non\\
&\quad\mbox{and strongly in }\,
C^0(\ov Q),\\
\label{conwt}
\dt\sialn\to\dt w^0&\quad\mbox{weakly-star in }\,L^\infty(Q).
\end{align}
\end{theorem}
\begin{proof}
By virtue of the global estimate \eqref{ssbound}, it follows the existence of a subsequence, which we label 
again by $n\in\enne$, and of limits $(\varphi,w,\xi)$ such that the convergence statements
\eqref{conphi}--\eqref{conwt} hold true with $(\vp^0,w^0,\xi^0)$ replaced by $(\vp,w,\xi)$. In this
connection, the strong convergence results in \eqref{conphi} and \eqref{conw} follow from standard
compactness results (cf., e.g., \cite[Sect.~8, Cor.~4]{Simon}). Observe that \an{the strong convergence in \eqref{conphi} along with} the Lipschitz continuity of $F_2'$ implies that
$\,F_2'(\phialn)\to F_2'(\varphi)\,$ strongly in $\,C^0(\overline Q)$.

We \an{then} need to show that $(\vp,w,\xi)=(\vp^0,w^0,\xi^0)$. To this end,
consider the time-integrated version of the system \eqref{wf:1}--\eqref{wf:2} with test functions
$v\in L^2(0,T;V)$ for $F_1=\faln$ and control $u=u_{\gamma_n}$ for $n\in \enne$. Passage to the
limit as $n\to\infty$ then shows that $(\vp,w,\xi)$ satisfies the initial conditions and 
is a solution to the time-integrated version
of \eqref{wf:1}--\eqref{wf:2} for the control $\,u$, which is equivalent to 
\eqref{wf:1}--\eqref{wf:2}. In order to conclude the proof, it remains to show that 
$\,\xi\in\sindi(\an{\ph})\,$ almost everywhere in $Q$. Indeed, if this is the case, then  $(\vp,w,\xi)$ is the
(uniquely determined) solution to the weak form of the state system for $F_1=\indi$ and
control $\,u\,$ and thus coincides with
$(\vp^0,w^0,\xi^0)$. Once this is shown, the unicity of the limit point also entails that
the convergence properties \eqref{conphi}--\eqref{conwt} are actually valid for the whole sequence
$\{\gamma_n\}$ and not just for a subsequence. 
 
Now define on $L^2(Q)$ the convex functional
\begin{equation*}
\Phi(v)=\int_Q \indi(v), \quad\mbox{if \,$\indi(v)\in L^1(Q),\,$ \,\,and 
\,$\Phi(v)=+\infty\,$,  \,otherwise}.
\end{equation*} 
It then suffices to show that $\xi$ belongs to the subdifferential of $\,\Phi\,$
at $\,\vp$, i.e., that
\begin{equation}
\label{reicht}
\Phi(v)-\Phi(\vp)\,\ge\,\int_Q \xi\,(v-\vp)\quad\forall\,v\in L^2(Q).
\end{equation}
At this point, we recall that $\phialn(x,t)\in [-1,1]$
\an{in} $\overline Q$.
Hence, by \eqref{conphi}, also $\vp(x,t)\in [-1,1]$ \an{in} $\overline Q$, and thus 
$\,\Phi(\vp)=0$. Now observe that in case that \,$I_{[-1,1]}(v)\not\in L^1(Q)$ \,the 
inequality \eqref{reicht} holds true since its left-hand side is infinite. If,
however,
$I_{[-1,1]}(v)\in L^1(Q)$, then obviously $v(x,t)\in [-1,1]$
almost everywhere in $\,Q$, and by virtue of \eqref{monoal} and \eqref{limF1al} 
it follows from Lebesgue's dominated convergence theorem that
$$
\lim_{n\to\infty}\int_Q \faln(v)= \Phi(v)=0.
$$
Now, by the convexity of $\faln$, and since $\faln(\phialn)$ is nonnegative, 
we have for all $v\in L^2(Q)$ that
$$
\faln'(\phialn)(v-\phialn)\,\le\,\faln(v)-\faln(\phialn)\,\le\,\faln(v)
\quad\mbox{a.e. in \,$Q$.}
$$
Using \eqref{conphi} and \eqref{conxi}, we thus obtain the following chain of (in)equalities:
\begin{eqnarray*}
\int_Q \xi(v-\vp) &\!\!=\!\!&\lim_{n\to\infty}\int_Q \faln'(\phialn)(v-\phialn)
\,\le\,{\limsup_{n\to\infty}}\int_Q\Bigl(\faln(v)-\faln(\phialn)\Bigr)\\
&\!\!\le\!\!&\lim_{n\to\infty}\int_Q\faln(v)\,=\,\Phi(v)\,=\,\Phi(v)-\Phi(\vp),
\end{eqnarray*}
which shows the validity of \eqref{reicht}. 
This concludes the proof of the 
assertion.
\end{proof}
\begin{remark} {\em
Note that the stronger conditions on the data 
required by \eqref{initials}
yield additional regularity 
for the solution also in the case $F_1=\indi$ with respect to the one 
obtained from Theorem~\ref{THM:EX:WEAK}. Indeed, we have
\begin{eqnarray*}
&&\varphi^0 \in {\W{1,\infty} H \cap \H1 V} \cap \L\infty{W},\quad
\xi^0\in L^\infty(0,T;H),\\
&&w^0\in {\H2 H \cap \W{1,\infty} V \cap \H1 {W}},\quad
\dt w^0\in L^\infty(Q).
\end{eqnarray*}}
\end{remark}

The following result provides a qualitative comparison between the solutions associated with
$F_1=\indi$ and  $F_1=F_{1,\gamma}$ for $\gamma\in (0,1]$.
\begin{theorem}
Suppose that \ref{ass:1:constants}--\ref{ass:3:initials} hold, and let $(\vp^0,w^0,\xi^0)=\SO(u)$ 
and, for any $\gamma\in(0,1]$, $(\vp_\gamma,w_\gamma, \pier{F'_{1,\gamma}(\vp_\gamma)})=\S_{\gamma}(u)$. Then it holds the estimate 
\begin{align}
	& \non
	\|\vp_{\gamma}- \vp^0\|_{\L\infty H \cap \L2 V}
	+ 	\|w_{\gamma}- w^0\|_{\H1 H \cap \L\infty V	}\,	\leq K_2 \,\gamma^{1/2}\,, 	
	\end{align}
with the positive constant $K_2$ introduced in Theorem~\ref{THM:CD:LOG}.
\end{theorem}
\begin{proof}
The result follows immediately from Theorem~\ref{THM:CD:LOG} and the semicontinuity properties of norms 
if we set $u_{\gamma_1}=u_{\gamma_2}=u$ and $\gamma_2=\gamma$
in the estimate \eqref{estiga12} and take the limit as $\gamma_1\searrow 0$.
\end{proof}

\section{Existence and approximation of optimal controls}

\setcounter{equation}{0}
Beginning with this section, we investigate the optimal control problem \CP0\ of minimizing the cost 
functional \eqref{cost} over the admissible set $\Uad$ subject to \an{the} state system 
\Sys\ in the form \eqref{wf:1},\eqref{wf:2},\eqref{ss4} for $F_1=\indi$ under the following
additional assumptions on the data of the cost functional:

\begin{enumerate}[label={\bf (A\arabic{*})}, ref={\bf (A\arabic{*})}, start=4]
\item \,
The constants $k_1,\ldots,k_6,\ell$  are nonnegative and not all equal to zero.
\label{ass:4}

\item \,
$\vp_\Omega, w_\Omega, w_\Omega'\in\Ldue$\, and \,$\vp_Q,w_Q,w_Q'\in L^2(Q)$.
\label{ass:5}

\item \,
$u_*,u^*\in L^\infty(Q)$ \an{and} satisfy $\,u_*\le u^*\,$ a.e. in $Q$. 
\label{ass:6}
\end{enumerate}

\vspace{1mm}\noindent 
In comparison with \CP0, 
we consider for $\gamma>0$ the following control problem:

\vspace{1mm}\noindent
\CPal \,\,\,Minimize $\,{\cal J}((\vp,w),u)\,$
for $\,u\in\Uad$, where  $\varphi,w$ \pier{denote the components}\\
\phantom{\CPal \,\,}
\pier{of the solution $\S_\gamma(u)$ to the state system.}

\vspace{1mm}\noindent
We expect that the minimizers  of \CPal\ are for $\gamma\searrow0$ related to minimizers of \CP0.
Prior to giving an affirmative answer to this conjecture, we first recall that \CPal\ has\pier{,} by virtue
of \cite[Thm.~3.1]{CSS3}\pier{,} for every $\gamma>0$ a solution; a corresponding result for \CP0\ is not
yet known and will be shown below. We  begin our analysis with the following result.

\begin{proposition}\label{PROP:J}
Suppose that \ref{ass:1:constants}--\ref{ass:6} and \eqref{initials} are satisfied, and let 
sequences $\,\{\gamma_n\}\subset (0,1]\,$ and
$\,\{u_n\}\subset\Uad\,$ be given such that, as $n\to\infty$, $\,\gamma_n\searrow0\,$ and 
$\,u_n\to u\,$ weakly-star in ${\cal U}$
for some $\,u\in\Uad$. Then \an{it can be shown that}
\begin{align}
\label{cesareuno}
&{\mathcal{J}}(\SO(u),u)\,\le\,\liminf_{n\to\infty}\,{\mathcal{J}}
({\cal S}_{\gamma_n}(u_n),u_n),\\[0.5mm]
\label{cesaredue}
&{\mathcal{J}}(\SO(v),v)\,=\,\lim_{n\to\infty}\,
{\mathcal{J}}({\cal S}_{\gamma_n}(v),v) \quad\forall\,v\in\Uad.
\end{align}
\end{proposition}

\begin{proof}
Theorem~\ref{THM:DQ} yields that $(\phialn,\sialn\pier{,F'_{1, \gamma_n}(\phialn)})=\S_{\gamma_n}(u_n)$  fulfills the convergence relations
\pier{\eqref{conphi}--\eqref{conw}}. 
The validity of \eqref{cesareuno} is then a direct consequence of 
the semicontinuity properties of the cost
functional \an{$\J$}.
 
Now suppose that $v\in\Uad$ is arbitrarily chosen, and put 
$(\phialn,\sialn\pier{,F'_{1, \gamma_n}(\phialn)}):=\S_{\gamma_n}(v)$ for all $n\in\enne$,
 as well as $(\vp^0,w^0,\xi^0):=\SO(v)$.
Applying Theorem~\ref{THM:DQ} with the constant sequence $u_n=v$, $n\in\enne$, we see that
\eqref{conphi}--\eqref{conwt} are valid once more. Since the first six summands of the cost functional are continuous with respect to the strong topology of $C^0([0,T];H)$, we
conclude the validity of
\eqref{cesaredue}. 
\end{proof}

We are now in a position to prove the existence of minimizers for the \an{control} problem \CP0. We have the following result.

\begin{corollary}\label{COR}
Suppose that \ref{ass:1:constants}--\ref{ass:6} and \eqref{initials} are fulfilled. Then the optimal control problem 
\CP0\ has at least
one solution.
\end{corollary}

\begin{proof}
Pick an arbitrary sequence $\{\gamma_n\}\subset (0,1]$ such that $\gamma_n\searrow0$ as $n\to\infty$.
Then the optimal control problem {\bf (CP$_{\boldsymbol{\gamma_n}}$)} has for every $n\in\enne$ a solution 
$((\phialn,\sialn),u_{\gamma_n})$, where $(\phialn,\sialn\pier{,F'_{1, \gamma_n}(\phialn)})=\S_{\gamma_n}(u_{\gamma_n})$ 
for $n\in\enne$.
Since $\Uad$ is bounded in ${\cal U}$, we may without loss of generality 
assume that $u_{\gamma_n}\to u$ weakly-star in ${\cal U}$
for some $u\in\Uad$\an{, the latter being a consequence of the convexity and the strong \juerg{closedness} of $\Uad$}.  We then obtain from Theorem~\ref{THM:DQ} that \eqref{conphi}--\eqref{conwt} 
  hold true with $(\vp^0,w^0,\xi^0)=\SO(u)$. Invoking the optimality of
$((\phialn,\sialn),u_{\gamma_n})$ for (${\mathcal{CP}}_{\gamma_n}$), we 
then find from Proposition~\ref{PROP:J} for every $\,\an{v}\in\Uad\,$ the chain of (in)equalities
\begin{align}
&{\cal J}(\SO(u),u)\,\le\,\liminf_{{n}\to\infty}\,
{\cal J}(\S_{\gamma_n}(u_{\gamma_n}),u_{\gamma_n})\,\le\,
\liminf_{{n}\to\infty}\,{\cal J}(\S_{\gamma_n}(\an{v}),\an{v})
\,=\,{\cal J}(\SO(\an{v}),\an{v} ),\non
\end{align}
which yields that $\,(\SO(u),u)\,$ is an optimal pair for \CP0. The assertion is thus proved.
\end{proof}

Theorem~\ref{THM:DQ} and the proof of Corollary~\ref{COR} indicate that optimal controls of \CPal\ 
are ``close'' to optimal
controls of \CP0\ as $\gamma$ approaches zero. However, they do not yield any information on 
whether every optimal control of \CP0\ can be approximated in this way. In fact, such a global result cannot be expected to hold true. 
Nevertheless, a local answer can be given by employing a well-known trick. To this end, 
let $\ubar\in\Uad$ be an optimal control for \CP0\ with the associated state 
$\SO(\ubar)$. We associate
with this optimal control the {\em adapted cost functional}
\begin{equation}
\label{adcost}
\widetilde{\cal J}((\vp,w),u):=
{\cal J}((\vp,w),u)\,+\,\frac 12\,\|u-\ubar\|^2_{\an{L^2(Q)}}
\end{equation}
and a corresponding \emph{adapted optimal control problem} for $\gamma>0$, namely:

\vspace{2mm}\noindent
\CPaltil\quad Minimize $\,\, 
\widetilde {\cal J}((\vp,w),u)\,\,$
for $\,u\in\Uad$\,  subject to \,$(\vp,w)=\SAL(u)$.

\vspace{2mm}\noindent
With essentially the same proof as that of \cite[Thm.~3.1]{CSS3} (which needs no repetition here), we can show that the 
adapted optimal control problem 
${\bf (\widetilde{CP}_\gamma)}$ has for every $\gamma>0$ at least one solution.
The following result gives a partial answer to the question raised above
\an{concerning the approximation of optimal controls for \CP0\ by the approximating problem \CPaltil}.

\begin{theorem} \label{THM:APPROX}
Let the assumptions of Proposition~\ref{PROP:J} be fulfilled, suppose that 
$\ubar\in \Uad$ is an arbitrary optimal control of \CP0\ with associated state  
$(\bvp,\overline w,\overline\xi)=\S_0(\ubar)$, and let $\,\{\gamma_k\}_{k\in\enne}\subset 
(0,1]\,$ be
any sequence such that $\,\gamma_k\searrow 0\,$ as $\,k\to\infty$. 
\pier{Then, for any $k\in \enne$ there exists 
an optimal control
 $\,u_{\gamma_{k}}\in \Uad\,$} of the adapted problem 
\pier{${\bf (\widetilde{CP}}_{\boldsymbol{\gamma_{k}}}{\bf )}$}
 with associated state \pier{$(\vp_{\gamma_k},w_{\gamma_k},\xi_{\gamma_k})=\S_{\gamma_{k}}
(u_{\gamma_{k}})$,
 such that, as $k\to\infty$,}
\begin{align}
\label{tr3.4}
&\pier{u_{\gamma_{k}}}\to \ubar \quad\mbox{strongly in }\,L^2(Q),
\end{align}
and such that \eqref{conphi}--\eqref{conwt} hold true with $(\vp^0,w^0,\xi^0)$ replaced 
by $(\bvp,\overline w,\overline\xi)$. 
 Moreover, we have 
\begin{align}
\label{tr3.5}
&\lim_{n\to\infty}\,\widetilde{{\cal J}}(\S_{\gamma_{n}}(u_{\gamma_{n}}),u_{\gamma_{n}})
\,=\,{\cal J}(\SO(\ubar),\ubar).
\end{align}
\end{theorem}

\begin{proof}
For any $ k\in\enne$, we pick an optimal control
$u_{\gamma_k} \in \Uad\,$ for the adapted problem ${\bf (\widetilde{CP}}_{\boldsymbol{\gamma_k}}{\bf )}$ and denote by 
$(\vp_{\gamma_k},w_{\gamma_k},\pier{\xi_{\gamma_k}})=\S_{\gamma_k}(u_{\gamma_k})$ the associated strong solution to the 
state system \Sys .
By the boundedness of $\Uad$ in $\calU$, there is some subsequence $\{\gamma_{n}\}$ 
of $\{\gamma_k\}$ such that
\begin{equation}
\label{ugam}
u_{\gamma_{n}}\to u\quad\mbox{weakly-star in }\,{\cal U}
\quad\mbox{as }\,n\to\infty,
\end{equation}
for some $u\in\Uad$. 
Thanks to Theorem~\ref{THM:DQ}, the convergence properties \eqref{conphi}--\eqref{conwt} hold true correspondingly
for the triple $(\vp^0,w^0,\xi^0)=\S_0(u)$. 
In addition,  the pair $(\S_0(u),u)$
is admissible for \CP0. 

We now aim at showing that $ u=\ubar$. Once this is shown, it follows from the unique solvability of
the state system that also $(\vp^0,w^0,\xi^0)=
(\bvp,\overline w, \overline\xi)$. 
Now observe that, owing to the weak sequential lower semicontinuity properties of 
$\widetilde {\cal J}$, 
and in view of the optimality property of $(\SO(\ubar),\ubar)$ 
for problem \CP0,
\begin{align}
\label{tr3.6}
\liminf_{n\to\infty}\, \widetilde{\cal J}(\S_{\gamma_n}(u_{\gamma_n}),u_{\gamma_n})\,
&\ge \,{\cal J}(\SO(u),u)\,+\,\frac{1}{2}\,
\|u-\ubar\|^2_{\an{L^2(Q)}}\nonumber\\[1mm]
&\geq \, {\cal J}(\SO(\ubar),\ubar)\,+\,\frac{1}{2}\,\|u-\ubar\|^2_{\an{L^2(Q)}}\,.
\end{align}
On the other hand, the optimality property of  
$\,(\S_{\gamma_{n}}(u_{\gamma_n}),u_{\gamma_n})
\,$ for problem ${\bf (\widetilde {CP}}_{\boldsymbol{\gamma_{n}}}{\bf )}$ yields that
for any $n\in\enne$ we have
\begin{equation}
\label{tr3.7}
\widetilde {\cal J}({\cal S}_{\gamma_{n}}(u_{\gamma_{n}}),
u_{\gamma_{n}})\,\le\,\widetilde {\cal J}({\cal S}_{\gamma_n}
(\ubar),\ubar)\,{{}=\, {\cal J}({\cal S}_{\gamma_n}
(\ubar),\ubar)\,,}
\end{equation}
whence, taking the limit superior as $n\to\infty$ on both sides and invoking (\ref{cesaredue}) in
Proposition~\ref{PROP:J},
\begin{align}
	\non
	&\limsup_{n\to\infty}\,\widetilde {\cal J}(\S_{\gamma_{n}}(u_{\gamma_n}),
u_{\gamma_n})\,\le\,\limsup_{n\to\infty}\widetilde {\cal J}(\S_{\gamma_n}(\ubar),\ubar) 
\\ & \quad 
\label{tr3.8}
=\,\limsup_{n\to\infty} {\cal J}(\S_{\gamma_n}(\ubar),\ubar)
\,=\,{\cal J}(\SO(\ubar),\ubar)\,.
\end{align}
Combining (\ref{tr3.6}) with (\ref{tr3.8}), we have thus shown that 
$\,\frac{1}{2}\,\|u-\ubar\|^2_{L^2(Q)}=0$\,,
so that $\,u=\ubar\,$  and thus also $(\bvp,\overline w,\overline\xi)
=(\vp^0,w^0,\xi^0)$. 
Moreover, (\ref{tr3.6}) and (\ref{tr3.8}) also imply that
\begin{align*}
&{\cal J}(\SO(\ubar),\ubar) \, =\,\widetilde{\cal J}(\SO(\ubar),\ubar)
\,=\,\liminf_{n\to\infty}\, \widetilde{\cal J}(\S_{\gamma_n}(u_{\gamma_n}),
 u_{\gamma_{n}})\nonumber\\[1mm]
&\,=\,\limsup_{n\to\infty}\,\widetilde{\cal J}(\S_{\gamma_n}(u_{\gamma_n}),
{u_{\gamma_n}}) \,
=\,\lim_{n\to\infty}\, \widetilde{\cal J}(\S_{\gamma_n}(u_{\gamma_n}),
u_{\gamma_n})\,,
\end{align*}                                    
which proves {(\ref{tr3.5})}. 
Moreover, the convergence properties 
\eqref{conphi}--\eqref{conwt} are satisfied. On the other hand, we have that  
\begin{align*}
{\cal J}(\SO(\ubar),\ubar) 
\,&\leq\,\liminf_{n\to\infty}\, {\cal J}(\S_{\gamma_n}(u_{\gamma_n}),
 u_{\gamma_{n}})
 \,\leq\,\limsup_{n\to\infty}\, {\cal J}(\S_{\gamma_n}(u_{\gamma_n}),
 u_{\gamma_{n}}) \nonumber\\[1mm]
&\leq \,\limsup_{n\to\infty}\,\widetilde{\cal J}(\S_{\gamma_n}(u_{\gamma_n}),
{u_{\gamma_n}}) \,
=\,{\cal J}(\SO(\ubar),\ubar) ,
\end{align*}        
so that also 
$ {\cal J}(\S_{\gamma_n}(u_{\gamma_n}),{u_{\gamma_n}})$ converges to $ {\cal J}(\SO(\ubar),\ubar)$ as 
$n\to \infty $, and the relation in \eqref{adcost} enables us to infer \an{the strong convergence in} \eqref{tr3.4}
\juerg{for the subsequence $\{u_{\gamma_n}\}$.}  

\juerg{We now claim that \eqref{tr3.4} holds true even for the entire sequence, due to the complete identification of the limit $u$ as $\ubar$. 
Assume that \eqref{tr3.4} does not hold true. Then there exist some $\varepsilon>0$ and a subsequence $\{\gamma_j\}$ of  $\{\gamma_k\}$ such that  
$$
\|u_{\gamma_j}-\ubar\|_{L^2(Q)}\,\ge\,\varepsilon \quad\forall\,j\in\enne.
$$
However, by the boundedness of $\Uad$, there is some subsequence $\{\gamma_{j_n}\}$ 
of $\{\gamma_j\}$ such that, with some $\tilde u\in\Uad$,
\begin{equation*}
u_{\gamma_{j_n}}\to \tilde u\quad\mbox{weakly-star in }\,{\cal U}
\quad\mbox{as }\,n\to\infty\,.
\end{equation*}
Arguing as above, it then turns out that $\tilde u = \ubar $ and that \eqref{tr3.4} holds for the subsequence
$\{u_{\gamma_{j_n}}\}$ as well, which contradicts the obvious fact that $\{ u_{\gamma_{j}} \}$ cannot have a subsequence which 
converges strongly to $\,\ubar\,$ in $L^2(Q)$.}
\end{proof}


\section{First-order \an{n}ecessary \an{o}ptimality \an{c}onditions}
\label{SEC:FOC}
\setcounter{equation}{0}

We now derive first-order necessary optimality conditions for
the control problem \CP0, using the corresponding conditions
for \CPaltil\ as approximations. To this end, we generally assume that the conditions
\ref{ass:1:constants}--\ref{ass:6} and \eqref{initials} are fulfilled. \an{Moreover, w}e need an additional assumption:

\begin{enumerate}[label={\bf (A\arabic{*})}, ref={\bf (A\arabic{*})}, start=7]
\item \label{ass:7}
\,At least one of the conditions $\,k_6=0\,$ or $\,w_\Omega'\in V$\, is satisfied.
\end{enumerate}
Now let $\ubar\in\Uad$ be any fixed optimal control for \CP0\ with
associated state $(\bvp,\overline w,\overline\xi)\linebreak =\SO(\ubar)$, and assume that $\gamma\in(0,1]$ \an{is} fixed. 
Moreover, \an{suppose} that $\ubar_\gamma\in\Uad$ is an optimal control for \CPaltil\ with corresponding state
$(\bvp_\gamma,\overline w_\gamma)=\SAL(\ubar_\gamma)$. The corresponding adjoint problem is given, in its strong form \an{for simplicity}, by
\begin{alignat}{2}
	\non
	& - \dt \pal - F_2'(\bvp_\gamma)\, \dt \qal 
	- \Delta \pal  + \fal''(\bvp_\gamma)\, \pal + 
	\tfrac 2 \thc  F_2''(\bvp_\gamma)\, \pal
	- \tfrac 1{\thc^2} \dt\bw_\gamma \, F_2''(\bvp_\gamma)\, \pal &&
 	\\ & \qquad \label{adsys1}
 	= k_1 (\bvp_\gamma - \vp_Q) \quad&&\mbox{in }\,Q,\\
	& \non
	-\dt \qal - \a \Delta \qal +  \b \Delta (1 \circledast \qal) - \tfrac 1{\thc^2} F_2'(\bvp_\gamma)\, \pal &&
	\\ & \qquad \label{adsys2}
	= k_3 ( 1 \circledast (\bw_\gamma - w_Q) )
	+ k_5 (\dt \bw_\gamma - w_Q')
	+ k_4 (\bw_\gamma(T) - w_\Omega)
	  \quad&&\mbox{in }\,Q,\\
\label{adsys3}
& \dn \pal= \dn \qal = 0\quad&&\mbox{on }\,\Sigma,\\[2mm]
\non
& \pal(T)= {k_2}(\bvp_\gamma(T) - \vp_\Omega) - {k_6}F_2'(\bvp_\gamma(T)) (\dt \bw_\gamma(T) - w'_\Omega) , &&
	\\ \label{adsys4}
 & \quad  \qal(T)={k_6} (\dt \bw_\gamma(T) - w'_\Omega) \quad &&\mbox{in }\,\Omega,
\end{alignat}
\Accorpa\Adj {adsys1} {adsys4}
where the product $\circledast $ is defined in \eqref{intr2}. Let us, for convenience, denote by $f_{\qal}$ the source term in \eqref{adsys2}, that is,
\begin{align}\label{source}
	f_{\qal}:= k_3 ( 1 \circledast (\bw_\gamma - w_Q) )
	+ k_5 (\dt \bw_\gamma - w_Q')
	+ k_4 (\bw_\gamma(T) - w_\Omega). 
\end{align}
According to \cite[Thm.~3.6]{CSS3}, the adjoint system has under the assumptions \ref{ass:1:constants}--\ref{ass:7} a unique weak 
solution
\begin{align}
	\label{reg:adj:p}
	\pal & \in \H1 \Vp \cap \L\infty H \cap \L2 V,
	\\
	\label{reg:adj:q}
	\qal & \in \H1 H \cap \L\infty V \cap \L2 W,
\end{align}
that satisfies the weak variational form 
\begin{align}
	& - \< \dt \pal, v>_V
	- \iO F_2'(\bvp_\gamma)\, \dt \qal \, v
	+ \iO \nabla \pal \cdot \nabla v  
	+ \iO \fal''(\bvp_\gamma) \,\pal \, v
		\non \\ & \quad 
	+ \frac 2 \thc \iO F_2''(\bvp_\gamma)\,\pal \, v
	- \frac 1{\thc^2} \iO \dt\bw_\gamma \,F_2''(\bvp_\gamma)\, \pal \, v
	= \an{k_1}\iO  (\bvp_\gamma - \vp_Q)\, v,
	\label{pier31}\\
	& - \iO \dt \qal \, v
	+ \a \iO \nabla  \qal \cdot \nabla v 
	- \b  \iO \nabla (1 \circledast \qal ) \cdot \nabla v
	- \frac 1{\thc^2}  \iO F_2'(\bvp_\gamma) \pal \, v
= \iO f_{\qal} \,v \label{pier32}
	\end{align}
for every $v\in V$, almost everywhere in $(0,T)$, and 
the final conditions
\begin{alignat}{2}
	\pal(T) & = {k_2}(\bvp_\gamma(T) - \vp_\Omega) - {k_6}F_2'(\bvp_\gamma(T)) (\dt \bw_\gamma(T) - w'_\Omega) \quad && \an{\text{in $\Omega$,}} \label{pier33}
		\\ 
 q_\gamma(T) & ={k_6} (\dt \bw_\gamma(T) - w'_\Omega)\label{pier34}
   &&\an{\text{in $\Omega$}}.
\end{alignat}

The variational inequality representing the first-order necessary optimality condition for \CPaltil\
then takes the form
\begin{equation}
\label{vuggamma}
\int_Q \bigl(\qal \,+\,\ell \,\ubar_\gamma \,+\,(\ubar_\gamma-\ubar)\bigr)(v-\ubar_\gamma)\,\ge\,0\quad\forall\,v\in\Uad. 
\end{equation}

In the following, we now derive some a priori bounds for the adjoint state variables, where we denote by $C_i$, $i\in\enne$,
constants that depend only on the data of the problem and not on $\gamma\in (0,1]$. At first, we conclude from the general
assumptions \ref{ass:1:constants}--\ref{ass:7} \an{and \eqref{initials}} and from the global bound
\eqref{ssbound} that the following holds true:
\begin{align}
\label{boundal1}
&\|F_2'(\bvp_\gamma)\|_{C^0(\ov Q)}\,+\,\|F_2''(\bvp_\gamma)\|_{C^0(\ov Q)}\,+\,
\|k_1 (\bvp_\gamma - \vp_Q)\|_{L^2(Q)}\,+\,
\|f_{\qal}\|_{L^2(Q)}\nonumber\\
&+\,\an{\| \pal(T)\|}\,+\,\| \qal(T)\|_{V}\,\le\,C_1\quad\forall\,\gamma\in(0,1]\,.
\end{align}

\noindent
{\bf First estimate:}
We {take $v=p\an{_\gamma}$ in \eqref{pier31}, $v = - \thc^2 \dt q\an{_\gamma}$ in \an{\eqref{pier32}},}
add the resulting equalities and note that two terms cancel out. Then, we integrate over $(t,T)$ and by parts. Putting
$\,Q_t^T:=\Omega\times (t,T)$, we obtain the identity
\begin{align}
	& \non
	 \frac 12\IO2 \pal
	+ \Qtt |\nabla \pal|^2
	+ \Qtt \fal''(\bvp_\gamma) |\pal|^2
	+ \thc^2 \Qtt |\dt \qal|^2	
	+\frac {\a \thc^2}{2} \IO2 {\nabla \qal}
	\\ & \quad \non
		= 
	\frac 12\IOT2 \pal
	+\frac {\a \thc^2}{2} \IOT2 {\nabla \qal}
	+k_1 \Qtt(\bvp_\gamma - \vp_Q) \pal
	- \frac 2 {\thc} \Qtt F_2''(\bvp_\gamma) \, \pal^2
	\\ & \qquad 
	+ \frac 1 {\thc^2} \Qtt\dt \bw_\gamma \,F_2''(\bvp_\gamma) \, \pal^2
	\an{-} \b \thc^2 \Qtt \nabla (1 \circledast \qal ) \cdot \nabla (\dt \qal)
	- \thc^2 \Qtt f_{q_\gamma} \, \dt \qal.
	\label{adj:est}
\end{align}
Notice that the third term on the \lhs\ is nonnegative since \pier{$\fal''\ge 0$.} 
As for the sixth term on the \rhs, we note that $(1 \circledast \qal )(T)=0 $
in $\Omega$, thus the Young and \Holder\ inequalities allow us to deduce that
\begin{align*}
		&  \an{-} \b \thc^2 \Qtt \nabla (1 \circledast \qal ) \cdot \nabla (\dt \qal)		\\ 
		& \quad \an{=}  \b \thc^2 \iO  \nabla (1 \circledast \qal ) (t)  \cdot \nabla \qal(t)
		\an{-} \b \thc^2 \Qtt |\nabla \qal|^2
\\ 
&\quad{}	\leq 	
\frac {\a \thc^2}{4} \IO2 {\nabla \qal}
			+ C_2 \Qtt  |\nabla \qal|^2 .	
\end{align*}
For third and last terms on the \rhs, we infer from \eqref{boundal1} and  Young's inequality that
\begin{align*}
	 k_1 \Qtt(\bvp_\gamma - \vp_Q)\pal
- \thc^2 \Qtt f_{q_\gamma} \, \dt \qal
	\leq 
	\frac {\thc^2}2 \Qtt |\dt \qal|^2
	+ C_3 \an{\Qtt (|p_\gamma|^2+1)} \,.
\end{align*}
Moreover, \eqref{boundal1} implies that the terms involving the terminal conditions are bounded by a constant $C_4>0$. Finally,
we invoke \eqref{boundal1} and the fact that  $\dt \bw_\gamma$ is uniformly bounded in $L^\infty(Q)$ \an{to deduce} the estimate
\begin{align*}
	- \frac 2 {\thc} \Qtt F_2''(\bvp_\gamma) \, \pal^2
	+ \frac 1 {\thc^2} \Qtt\dt \bw_\gamma \,F_2''(\bvp_\gamma) \, \pal^2
	\leq C_\an{5} \Qtt |\pal|^2.
\end{align*}
Collecting the above computations, and applying Gronwall's lemma, we infer that 
\begin{align}\label{boundal2}
	\norma{\pal}_{\L\infty H \cap \L2 V}
	+\norma{\qal}_{\H1 H \cap \L\infty V}
	\leq C_\an{6} \quad\forall\,\gamma\in(0,1].
\end{align}

\noindent
{\bf Second estimate:}
Next, we proceed with comparison in equation \eqref{adsys2} to deduce that
\begin{align*}
	\big\|\Delta \big(\a \qal \pier{{}-{}} \b (1 \circledast \qal) \big)\big\|_{\L2 H } \leq C_\an{7}\quad\forall\,\gamma\in(0,1].
\end{align*}
Then, setting $g_\gamma = \a \qal \pier{{}-{}} \b (1\circledast \qal)$, the elliptic regularity theory entails that $\norma{g_\gamma}_{\L2 W} \leq C_\an{8} .$ Hence, solving the equation $\a \qal \pier{{}-{}} \b (1\circledast \qal)=g_\gamma$ with respect to $1\circledast \qal$, 
we eventually obtain~that
\begin{align}\label{boundal3}
	\norma{1 \circledast \qal }_{\L2 W} + \norma{\qal }_{\L2 W} \leq C_\an{9}\quad\forall\,\gamma\in (0,1]. 
	\end{align}
\noindent

\noindent{\bf Third estimate:}
For the next estimate, we introduce the space 
\begin{equation}
\label{defQ}
{\cal Q}=\{v\in H^1(0,T;V^*)\cap L^2(0,T;V): \ v(0)=0\}, 
\end{equation}
which is a closed subspace of $H^1(0,T;V^*)\cap L^2(0,T;V)$ and thus a Hilbert space. 
As is well known, ${\cal Q}$ is continuously embedded in $C^0([0,T];H)$, and we
have the dense and continuous embeddings ${\cal Q}\subset L^2(0,T;H)\subset {\cal Q}^*$, where it is understood that
\begin{equation}\label{embo}
\langle v,w\rangle_{\cal Q}\,=\,\int_0^T(v(t),w(t))\,dt \quad\mbox{for all \,$w\in{\cal Q}$\, and $\,v\in L^2(0,T;H)$}.
\end{equation}  
Next, we recall the well-known  integration-by-parts formula for functions in $H^1(0,T;V^*)\linebreak\cap L^2(0,T;V)$, which 
yields that for all $v\in \pier{{}{\cal Q}}$ it holds the estimate
\begin{align}
\label{adj8}
&\Big|\int_0^T\langle \dt\pal(t),v(t)\rangle_V\,dt \, \Big| 
\,\le\,\Big|\int_0^T \langle \dt v(t),\pal(t)\rangle_V\,dt\Big|\,+\,\Big|\iO \pal(T)\,v(T) \,-\,\iO \pal(0)\,v(0)\Big|
\nonumber\\[1mm]
&\quad \le\,\|\pal\|_{L^2(0,T;V)}\,\|\dt v\|_{L^2(0,T;V^*)}\,+\, \an{\|\pal(T)\|\,\|v(T)\|}
\nonumber\\[1mm]
&\quad \le\,{C_\an{10}}\,\|v\|_{H^1(0,T;V^*)}\,+\,{C_{\an{11}}}\,\|v\|_{C^0([0,T];H)}\,\le\,{C_{\an{12}}}\,\|v\|_{\cal Q},
\end{align}
where we used \eqref{boundal1} and \eqref{boundal2}. \an{T}his \an{actually} means that
\begin{equation}\label{boundal6}
\|\dt\pal\|_{{\cal Q}^*}\,\le\,C_{\an{13}} \quad\forall\,\gamma\in(0,1].
\end{equation}
At this point, we can conclude from the estimates \eqref{boundal1}, \eqref{boundal2}, \eqref{boundal3}, and \eqref{boundal6}, using 
\pier{a comparison argument} in \eqref{pier31}, that the linear \pier{and continuous} mapping 
\begin{equation}\label{defLamal}
\pier{\Lambda_\gamma:{\cal Q}  \to \erre, \quad
\langle \Lambda_\gamma,v\rangle_{\cal Q}: =
\int_Q \fal''(\bvp_\gamma)\,\pal\,v}\,,
\end{equation}
satisfies
\begin{equation}\label{boundal7}
\|\Lambda_\gamma\|_{{\cal Q}^*}\,\le\,{C_{\an{14}}}\quad\forall\,\gamma\in(0,1].
\end{equation}

\an{Consider now} any sequence $\gamma_n\searrow0$. According to Theorem~\ref{THM:DQ} and Theorem~\ref{THM:APPROX}, we may without loss of generality
assume that the sequence $\,\{\bu_{\gamma_n}\}\,$ converges strongly in $\,L^2(Q)\,$ to $\,\bu\,$ and that the convergence properties 
\eqref{conphi}--\eqref{conwt} are satisfied with $\,(\vp_{\gamma_n},w_{\gamma_n})\,$ and $\,(\vp^0,w^0,\xi^0)\,$ replaced by 
$\,(\bvp_{\gamma_n},\ov w_{\gamma_n})\,$ and $\,(\bvp,\ov w,\ov \xi)\,$, respectively. By virtue of the estimates
\eqref{boundal2}, \eqref{boundal3}, \eqref{boundal6}, and \eqref{boundal7}, we may also assume without loss of generality that\an{, as $n\to\infty$}
\begin{align}
\label{conp}
\paln\to p&\quad\mbox{weakly-star in $\,L^\infty(0,T;H)\cap L^2(0,T;V)$},\\
\label{conq}
\qaln\to q&\quad\mbox{weakly-star in }\,H^1(0,T;H)\cap L^\infty(0,T;V)\cap L^2(0,T;W),\\
\label{coniq}
1 \circledast \qaln\to 1\circledast q&\quad\mbox{weakly in }\,L^2(0,T;W),\\
\label{conpt}
\dt\paln \to\dt p&\quad\mbox{weakly in }\,{\cal Q}^*,\\
\label{conlam}
\Lambda_{\gamma_n}\to \Lambda&\quad\mbox{weakly in }\,{\cal Q}^*,
\end{align}
for suitable limit points \an{$p,q,$ and $\Lambda$}.

\an{Then we} perform a passage to the limit as $n\to\infty$ in the adjoint system \eqref{pier31}--\eqref{pier34},
written for $\gamma=\gamma_n$ and $\pier{(p,q)=(\paln,\qaln)}$, for $n\in\enne$. At first, we recall that by \eqref{conphi}
we have that $\,\bvp_{\gamma_n}\to \bvp\,$ strongly in $\,C^0(\ov Q)$, and \ref{ass:2:F2} implies that\an{, as $n\to\infty$}
\begin{equation}
\label{sigi}
F_2'(\bvp_{\gamma_n})\to F_2'(\bvp) \quad\mbox{and} \quad F_2''(\bvp_{\gamma_n})\to F_2''(\bvp), \quad\mbox{both strongly in }\,
C^0(\ov Q).
\end{equation}
From the convergence results
stated above it is then readily seen that, as $n\to\infty$,
\begin{alignat}{2}
\label{term1}
F_2'(\bvp_{\gamma_n})\,\dt\qaln & \to F_2'(\bvp)\,\dt q \quad && \mbox{weakly in }\,L^2(Q),\\
\label{term2}
F_2''(\bvp_{\gamma_n})\,\paln &\to F_2''(\bvp)\,p \quad && \mbox{weakly in }\,L^2(Q),\\  
\label{term3}
F_2'(\bvp_{\gamma_n})\,\paln &\to F_2'(\bvp)\,p  \quad && \mbox{weakly in }\,L^2(Q).
\end{alignat}
Next, observe that by virtue of \eqref{conw} \an{and \eqref{conwt}} we have that $\,\dt \ov w_{\gamma_n}\to\dt \ov w\,$ weakly-star in $H^1(0,T;H)\cap 
L^\infty(0,T;V)\cap L^2(0,T;W)\an{\cap L^\infty(Q)}$ and thus, by continuous \an{embedding}, also
\begin{align}\label{term4}
\an{\dt\ov w_{\gamma_n}\to \dt\ov w\quad\mbox{weakly in }\,C^0([0,T];V).}
\end{align}
In addition, \cite[Sect.~8, Cor.~4]{Simon} implies that we also may assume that
\begin{equation}
\label{vabene}
\dt \bw_{\gamma_n}\to \dt\bw \quad\mbox{strongly in }\,C^0([0,T];L^{\an{\sigma}}(\Omega))\quad\mbox{for $1\le \an{\sigma} <6$}.        
\end{equation}
It then easily follows from \eqref{source}, \eqref{pier33}, and \eqref{pier34}, that
\begin{align}\label{soso}
f_{q_{\gamma_n}}&\to f_q:=k_3 ( 1 \circledast (\bw - w_Q) )
	+ k_5 (\dt \bw - w_Q')
	+ k_4 (\bw(T) - w_\Omega)  \quad\mbox{weakly in }\,L^2(Q),\\
p_{\gamma_n}(T)&\to p(T)= k_2(\bvp(T) - \vp_\Omega) - k_6 F_2'(\bvp(T)) (\dt \bw(T) - w'_\Omega)  \quad\mbox{weakly in }\,L^2(\Omega),\\
q_{\gamma_n}(T)&\to q(T)=k_6 (\dt \bw(T) - w'_\Omega)  \quad\mbox{weakly in }\,V.
\end{align}
Finally, we claim that also
\begin{equation}\label{term6}
\dt\bw_{\gamma_n} \, F_2''(\bvp_{\gamma_n})\, p_{\gamma_n}\to \dt\bw \, F_2''(\bvp)\, p\quad\mbox{weakly in }\,L^2(Q).
\end{equation}
Indeed, we have for every $v\in L^2(Q)$ the identity
\begin{align*}
&\int_Q \bigl(\dt\bw_{\gamma_n}\,F_2''(\bvp_{\gamma_n})\,p_{\gamma_n}\,-\,\dt\bw \,F_2''(\bvp)\,p\bigr)\,v\\
&=\int_Q\bigl(\dt\bw_{\gamma_n}-\dt\bw\bigr)\,F_2''(\bvp)\,p_{\gamma_n}\,v\,
+\int_Q \dt\bw_{\gamma_n}\,\bigl(F_2''(\bvp_{\gamma_n})-F_2''(\bvp)\bigr)\,p_{\gamma_n}\,v\\
&\quad +\int_Q \dt\bw\,F_2''(\bvp)\,(p_{\gamma_n}-p)\,v\,=: I_{1n}+I_{2n}+I_{3n}\,,
\end {align*}
with obvious notation. \an{Since} $\dt\bw\,F_2''(\bvp)\,v\in L^2(Q)$, we have \pier{that} $I_{3n}\to 0$ as $n\to\infty$.
Moreover, the sequence $\{\dt\bw_{\gamma_n}\,p_{\gamma_n}\,v\}$ is bounded in $L^1(Q)$ so that \eqref{sigi} implies that
$I_{2n}\to 0$ as $n\to\infty$. Finally, using \eqref{boundal1}, \eqref{boundal2}, \eqref{vabene} \an{with $\sigma=4$}, H\"older's inequality, 
and the continuous embedding $V\subset L^4(\Omega)$, we see that
\begin{align*}
|I_{1n}|&\le \, C_{\an{15}}\int_0^T \|p_{\gamma_n}(t)\|_{\an{4}}\,\|\dt\bw_{\gamma_n}(t)-\dt\bw(t)\|_{\an{4}}\,\|v(t)\|\,dt\\ 
&\le\,C_{\an{16}}\,\|\dt\bw_{\gamma_n}-\dt\bw\|_{C^0([0,T];L^4(\Omega))}\,\|p_{\gamma_n}\|_{L^2(0,T;V)}\,\|v\|_{L^2(0,T;H)}
\,\,\,\to \,0\quad\mbox{as }\,n\to\infty,
\end{align*}
which proves the validity of the claim \eqref{term6}.

\an{Besides}, for every $v\in {\cal Q}$, 
\begin{align}
\label{susi}
\langle \dt p,v\rangle_{\cal Q}&\,=\,\lim_{n\to\infty}\langle \dt p_{\gamma_n},v\rangle_{\cal Q}
 \,=\,\lim_{n\to\infty} \int_0^T \langle \dt p_{\gamma_n}(t),v(t)\rangle _V\,dt \non\\
&\,=\,\lim_{n\to\infty}\Big( \iO p_{\gamma_n}(T)v(T)\, - \int_0^T\langle \dt v(t), p_{\gamma_n}(t)\rangle_V\,dt\Big)\non\\
&\,=\,\iO p(T)v(T) \,-\int_Q p\,\dt v\,.
\end{align}

At this point, we may pass to the limit as $n\to\infty$ in the adjoint system \eqref{pier31}--\eqref{pier34} to arrive at the following limit system:
\begin{align}
\label{wadj1neu}
&\langle \Lambda,v\rangle_{{\cal Q}}\,=\,-\int_Q p\,\dt v\,+\iO p(T) \,v(T) \,+\int_Q F_2'(\bvp)\,\dt q\,v
	- \int_Q\nabla p \cdot \nabla v
	\non\\
	&\quad - \frac 2{\theta_c} \int_Q F_2''(\bvp)\,p\,v\,+\,\frac 1{\theta_c^2}\int_Q\dt \bw \,F_2''(\bvp)\,p\,v
	\,+\int_Q k_1(\bvp-\vp_Q) \quad \an{\text{\and{for all} $\,v\in {\cal Q}$}}\,,
		\\[2mm]
	\label{wadj2neu}
	& -\iO \dt q(t)\,v \,+\,\alpha \iO \nabla q(t) \cdot \nabla v\,-\,\iO \nabla(1\circledast q(t))\cdot \nabla v\non\\
	&\quad -\frac 1{\theta_c^2}\iO F_2'(\bvp(t))\,p(t)\,v \,=\,\iO f_q(t)\,v
	\qquad\mbox{for all $\,v\in V$ and a.e. $t\in(0,T$)},
	\\[2mm]
	\label{wadj3neu}
	& 
	p(T)\,=\,k_2(\bvp(T)-\vp_\Omega)-k_6F_2'(\bvp(T))(\dt\bw(T)-w_\Omega')\quad\mbox{\an{ in }}\,\Omega,\\[2mm]
	\label{wadj4neu}
	&q(T)\,=\,k_6(\dt\bw(T)-w_\Omega') \quad\mbox{\an{in }}\,\Omega,
	\end{align}
where $f_q$ is defined in \eqref{soso}.

Finally, we consider the variational inequality \eqref{vuggamma} for $\gamma=\gamma_n$, $n\in\enne$. Passage to the
limit as $n\to\infty$, using the 
above convergence results, yields that
\begin{align}\label{gauss1}
\int_Q (q+\ell \,\ubar)(v-\ubar)\,\ge\,0 \quad\forall\,v\in \Uad.
\end{align}

Summarizing the above considerations, we have proved the following first-order necessary optimality conditions for the
optimal control problem \CP0.

\begin{theorem} 
Suppose that the conditions \ref{ass:1:constants}--\ref{ass:7} and \eqref{initials} are fulfilled, and let $\ubar\in\Uad$ be a 
minimizer of the optimal control problem \CP0\
with associate state $(\bvp,\bw,\ov \xi)=\SO(\ubar)$. Then there exist \an{$p,q$, and $\Lambda$}
such that the following holds true:\\[1mm]
{\rm (i)} \quad\,\,$p\in L^\infty(0,T;H)\cap L^2(0,T;V)$, $q\in H^1(0,T;H)\cap L^\infty(0,T;V)\cap L^2(0,T;W)$, 
$\Lambda\in {\cal Q}^*$.\\[1mm]
{\rm (ii)} \quad The adjoint system \eqref{wadj1neu}--\eqref{wadj4neu} and the variational inequality \eqref{gauss1}
are satisfied. 
\end{theorem}

\vspace{2mm}
\begin{remark}{\em
{(i)} Observe that the adjoint state $(p,q)$ and the Lagrange multiplier $\Lambda$ are not unique. However, all possible
choices satisfy \eqref{gauss1}, i.e., $\bu$ is for $\ell>0$ the $\an{L^2(Q)}$-orthogonal projection of $\an{-{\ell}^{-1}} q$ onto the 
closed and convex set $\{ u \in L^\infty(Q) : u_* \leq u \leq u^* \,\,\aeQ\}$, and 
\begin{align}
	\non
	\bu(x,t)=\max \big\{ 
	u_*(x,t), \min\{u^*(x,t),-{\an{\ell}}^{-1} q(x,t)\} 
	\big\} \quad \hbox{for $a.a. \,(x,t) \in Q.$}
\end{align} 
{(ii)} We have, for every $n\in\enne$, the complementarity slackness condition
{(cf.~\eqref{defLamal})}
$${\Lambda_{\gamma_n}(\paln)}=\int_Q F_{1,\gamma_n}''(\bvp_{\gamma_n})\,|\paln|^2
\pier{{}=\int_Q \frac{2\gamma_n}{1 - \bvp_{\gamma_n}^{\,2}}\,|\paln|^2}
\,\ge\,0.
$$
Unfortunately, \pier{our} convergence properties \pier{for  $\{\bvp_{\gamma_n}\}$
and} $\{\paln\}$ do not permit a passage to the limit in this
inequality to
derive a corresponding result for \CP0.}
\end{remark}

\section*{Acknowledgments}
This research was supported by the Italian Ministry of Education, 
University and Research~(MIUR): Dipartimenti di Eccellenza Program (2018--2022) 
-- Dept.~of Mathematics ``F.~Casorati'', University of Pavia. 
In addition, \pier{PC and GG aim to point out their collaboration,
as Research Associates, to the IMATI -- C.N.R. Pavia, Italy. Moreover, 
PC and AS gratefully acknowledge some support 
from the MIUR-PRIN Grant 2020F3NCPX ``Mathematics for industry 4.0 (Math4I4)'' 
and underline their affiliation}
to the GNAMPA (Gruppo Nazionale per l'Analisi Matematica, 
la Probabilit\`a e le loro Applicazioni) of INdAM (Isti\-tuto 
Nazionale di Alta Matematica).

\End{document}
